\documentclass[final]{siamltex}

\usepackage{amsmath}
\usepackage{amssymb,amsfonts,bm,bbm}
\usepackage{color,subfig}
\usepackage{xspace}
\usepackage{placeins}
\usepackage{rotating}
\usepackage{enumitem}

\newtheorem{thm}{Theorem}[section]

\newtheorem{remark}[thm]{Remark}

\newcommand{\N}{\mathbb N}
\newcommand{\R}{\mathbb R}

\newcommand{\bal}{{\bm \alpha}}
\newcommand{\bbe}{{\bm \beta}}
\newcommand{\bga}{{\bm \gamma}}

\newcommand{\bx}{{\bm x}}
\newcommand{\bz}{{\bm 0}}

\newcommand{\rM}{{\rm M}}
\newcommand{\rR}{{\rm R}}

\newcommand{\VnK}{V_{\!{\rm lo},n}(K)}
\newcommand{\VpK}{V_{\!{\rm ho},p}(K)}

\newcommand{\hVpK}{\widehat{V}_{\!{\rm ho},p}(K)}
\newcommand{\ProjnK}{\pi_{\!{\rm lo},n}^K}
\newcommand{\ProjpK}{\pi_{\!{\rm ho},p}^K}
\DeclareMathOperator*{\argmin}{arg\,min}

\newcommand{\bbmOne}{\mathbbm{1}}

\definecolor{shadecolor}{rgb}{0.9,0.9,0.9}

\newcommand{\pder}[2]{\ensuremath{\frac{\partial #1}{\partial #2}}} 

\title{A Discontinuous Galerkin method for Shock Capturing using
a mixed high-order and sub-grid low-order approximation space}

\author{Per-Olof Persson, Benjamin Stamm}

\begin{document}

\maketitle

\pagestyle{myheadings}
\thispagestyle{plain} 
\markboth{}{}

\begin{abstract}
This article considers a new discretization scheme for conservation laws. The discretization setting is based on a discontinuous Galerkin scheme in combination with an approximation space that contains high-order polynomial modes as well as piece-wise constant modes on a sub-grid.
The high-order modes can continuously be suppressed with a penalty function that is based on a sensor which is intertwined with the approximation space.
Numerical tests finally illustrate the performance of this scheme.
\end{abstract}

\begin{keywords} 
	Conservation laws,  
  	Discontinuous Galerkin Method,
  	High Order Methods,
  	Stability
\end{keywords}

\begin{abstract}
This article considers a new discretization scheme for conservation laws. The discretization setting is based on a discontinuous Galerkin scheme in combination with an approximation space that contains high-order polynomial modes as well as piece-wise constant modes on a sub-grid.
The high-order modes can continuously be suppressed with a penalty function that is based on a sensor which is intertwined with the approximation space.
Numerical tests finally illustrate the performance of this scheme.
\end{abstract}

\section{Introduction}

While discontinuous Galerkin (DG) and related high-order methods
\cite{cockburn01rkdg} are getting sufficiently mature to handle realistic
problems in many applications, there still exist challenging problems where state-of-the-art methods suffer from the lack of nonlinear stability and their high sensitivity to under-resolved features. 
This directly affects the solution of important problems involving shocks and turbulence models, but it
has also turned out to be a problem for simpler problems such as laminar or
inviscid flows, if the meshes are not well adapted to the solution
fields. This lack of robustness is one of the main challenges remaining for
the wide adoption of high-order methods.

Several approaches have been proposed to address the issue. One simple method
is to calculate a sensor that identifies the elements in the shock
region and reduce the degree of approximating polynomials
\cite{BauOden,burbeau01limiter}. This is usually combined with $h$-adaptivity
to better resolve the shocks \cite{dervieux03adaptation}, and it can be quite
satisfactory in particular for steady-state problems. 
In recent refined versions, these ideas have been combined with sub-cell resolution schemes, either by a high-order/low-order duality on unstructured~\cite{huerta2012shockDG} and structured~\cite{sonntag2017efficient} grids, or a finite-volume a-posteriori correction of marked cells on structured~\cite{dumbser2014subcell,vilar2019posteriori} and unstructured~\cite{dumbser2016simple} grids.
Another sophisticated class of approaches include limiting, for example based on weighted
essentially non-oscillatory (WENO) concepts \cite{eno1,weno1,weno2}. A related approach is to filter the solution in order to stabilize the discontinuities, see e.g. \cite{gassner2019shocks}. While these various
schemes have been demonstrated to handle very strong shocks with no robustness
issues, for various reasons they have not been widely employed to
large-scale 3D problems on unstructured meshes, which typically require
implicit solvers.

Alternatively, in \cite{persson06shock} it was demonstrated how a strategy
inspired by the early artificial viscosity methods can be very effective in
the context of high-order DG methods.  Based on the early stabilization
proposed by Neumann \cite{neumann50shocks} and used for discontinuous Galerkin
methods in an element-wise way
\cite{bassirebay95shocks,BauOden,hartmann02shock}, the method combines a
highly selective spectral sensor, based on orthogonal polynomials, with a
consistently discretized artificial viscosity added to the equations. The goal
is to smooth the discontinuities in the solution to a width that is
appropriately resolved by the mesh and the polynomial approximations, which
means in particular that the method obtains sub-cell resolution for high-order
discretizations, which gives a number of important benefits. The continuous
nature of the scheme allowed for Newton methods to produce fully converged
solutions to steady-state problems or to implicit time-stepping problems.  The
method has been widely employed and improved, e.g. in
\cite{klockner2011shock,barter08,persson13shock}.  A related class of schemes
uses so-called physics-based sensors and viscosities
\cite{jameson81jst,lele2009viscosity,moro14thesis}. Although popular, the
artificial viscosity approaches still suffer from spurious oscillations and
parameter sensitivities.

In the present work, we are proposing a new approach, based on a combination
of the excellent shock-capturing properties of the finite volume method and
the high-order accuracy of the DG method.  We define an approximation space
that contains both, high-order polynomial modes and low-order piece-wise
constant functions on a sub-grid, to design an approximation space.  The
fundamental idea is that the discontinuous Galerkin scheme, being a
projection-based scheme, can use both modes for defining the approximation and
the polynomial modes can be continuously suppressed by penalty if needed.
Appropriate numerical dissipation is automatically introduced by the jumps of
the solution.  The penalty, which is based on a sensor function is a delicate
issue.  However, the richness of the approximation space allows for a powerful
construction of such a sensor function which not only allows the detection of
shocks but also indicates efficiently recovery in an element where
penalization should be removed. 
Our approach is related to the contributions \cite{huerta2012shockDG,sonntag2017efficient,dumbser2014subcell,vilar2019posteriori,dumbser2016simple} referenced above involving a sub-grid as well. 
With exception of \cite{huerta2012shockDG}, these methods use a on/off-procedure to enable the (sub-grid) lower-order scheme based on an a posteriori approach for detection and reconstruction.
While the method in \cite{huerta2012shockDG} switches continuously between a low-order/high-order interpretation of the degrees of freedom we propose a combined approximation space and continuously suppress the high-order modes with a penalty function that arises naturally from the approximation space.

The outline of this article is as follows.  The second section introduces the
approximation space as well as several $L^2$-projections used within this
framework.  Section three considers space and time discretization of
(generally non-linear) conservation laws and Section four presents various
numerical results and tests in one and two spatial dimensions.
Section five is devoted to some conclusions.

\section{The approximation space and its properties}
Here, we first introduce the approximation space and propose some basis functions considered within this framework. Further, we will introduce different projections that will be used in the following.

\subsection{Polynomials with sub-grid components}

Our approach is a standard discontinuous Galerkin method with a particular
choice of approximation space. 
Consider a partition $T=\{ K_1, \ldots, K_{N_T} \}$ of a domain $\Omega\subset \mathbb R^d$, $d\ge 1$. 
Let $h_K$ denote the diameter of the element $K\in
T$.  In each element $K_\ell$, we define a sub-grid $T_{K_\ell}=\{k^\ell_1,\ldots,k^\ell_n\}$ of size $n$. 
Note that $n$ can in principle vary with $\ell$, this is however omitted for sake of a simple presentation. 
For a 1D
line segment, this sub-grid could for example be a uniform partition of each
element.  
Then, for each element $K\in T$, we first define the space of polynomials
of at most degree $p$ with zero average:
\begin{align}
	\VpK = \left \{ v\in \mathbb{P}_{\! p}(K) : \int_K v\,dx = 0 \right \},
\end{align}
where $\mathbb{P}_{\! p}(K)$ denotes a space of polynomials of total degree $p$ on
$K$, as well as the space of piece-wise constant functions in each element $k$
on the sub-grid triangulation $T_K$:
\begin{align}
\VnK = 
      \{ v\in L_2(K) : v|_k \in \mathbb{P}_{\! 0}(k)\  \forall k \in T_K \}.
\end{align}
Note that again, the polynomial degree $p$ can also vary locally in the partition $T$.
We then combine these spaces to obtain the approximation space that
 will be considered in this framework:
\begin{align}
  V_\delta(K) = \VpK \oplus \VnK.
\end{align}
Here, the two discretization parameters $p$ and $n$ are compactly denoted as $\delta:=(p,n)$. 
We note that the condition on the integral of the functions in the
polynomial space is necessary to avoid duplicating the constant
function on $K$.  Clearly, the space $V_\delta(K)$ reproduces the DG and the FVM
methods in the extreme cases:
\begin{itemize}
\setlength{\itemsep}{0pt}
\item If the sub-grid equals the original element ($n=1$ and thus $T_K = \{K\}$), the space
  is simply the polynomials of degree $p$ and the resulting space is the
  standard DG space.
\item If the polynomial degree is $p=0$, the space consists of the piece-wise
  constant functions on the sub-grid and the resulting space is the standard space that is used for the 
  FVM method.
\end{itemize}
However, for a non-trivial combination of $n$ and $p$, we obtain a new
space with interesting properties. Discontinuities can be handled efficiently by the piece-wise
constants, while the asymptotic high-order convergence rates for
smooth solutions are obtained by the polynomials.

 \begin{figure}[t]
   \begin{center}
     \includegraphics[width=.9\textwidth]{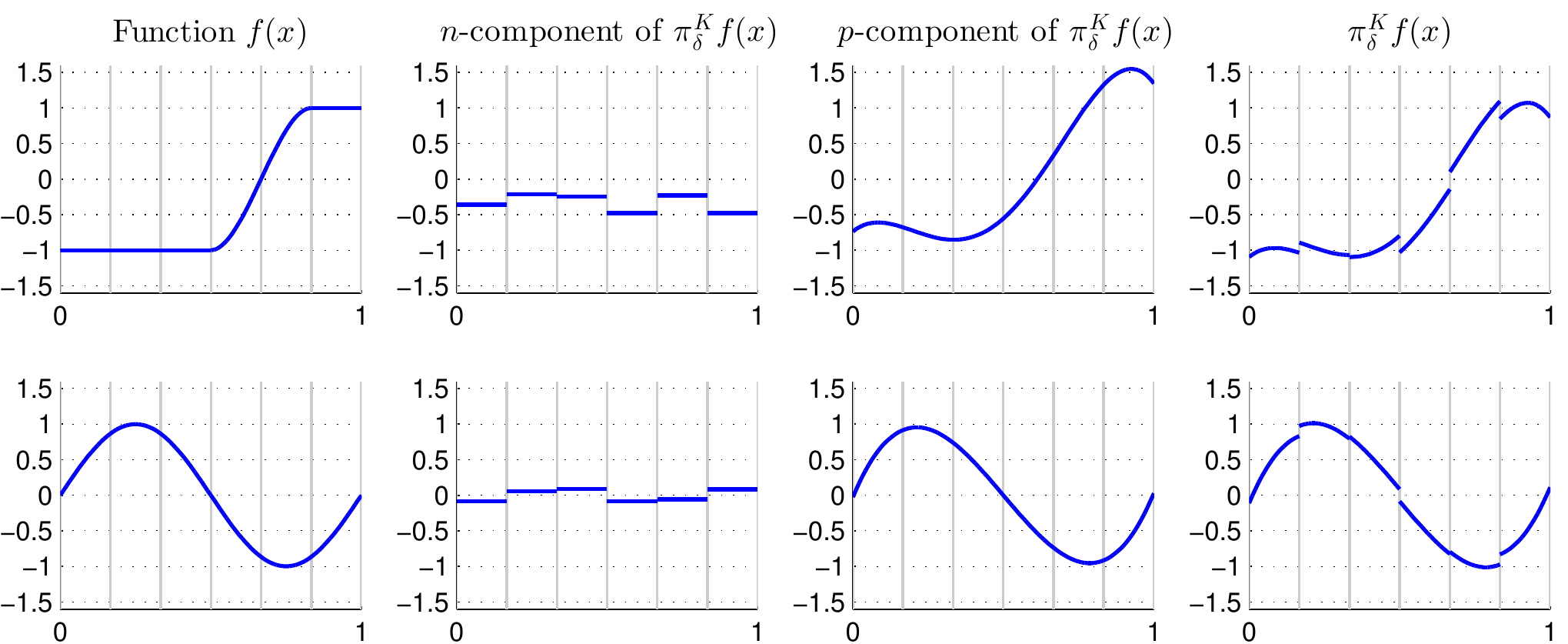}
   \end{center}
   \caption{Projection $\ProjpK f$ of two smooth example functions $f$ on the space $V_\delta(K)$ for $p=4$ and $n=6$ where the two components in $\VnK$ and $\VpK$ of $\ProjpK f$ are also indicated.}
   \label{fig:projection}
 \end{figure}

Finally, we now consider the global space given by
\[
	V_\delta = \bigoplus_{K \in T} V_\delta(K).
\] 
Therefore, any function $v_\delta \in V_\delta$ lies locally in $V_\delta(K)$, i.e. $v_\delta|_K \in V_\delta(K)$ for all $K\in T$.

\subsection{Basis functions}

Within each element $K\in T$, we introduce the following basis functions for
the space $V_\delta(K)$. Let $L^\ell_i$, $i=0,\ldots,N(p)$ be a set of orthogonal polynomials up to total degree $p$ on $K_\ell$ such that $L^\ell_0=1$. 
In one spatial dimension, we use the Legendre polynomials.
Then, we define the following $N(p)+n$ basis functions for $V_\delta(K)$ by
\begin{align}
  \varphi^\ell_i =
  \begin{cases}
    L^\ell_i & i=1,\ldots,N(p), \\
    \bbmOne_{k^\ell_j} & j = 1,\ldots,n,\ i=N(p)+j,
  \end{cases}
\end{align}
where we use the characteristic function
\begin{align}
  \bbmOne_{k}(x) =
  \begin{cases}
    1 & x \in k, \\
    0 & \text{otherwise}.
  \end{cases}
\end{align}
Note that while both, the
polynomial components and the low order piece-wise constant components, are orthogonal w.r.t the $L_2$-norm, the combined basis
is in general not orthogonal.

\subsection{Projections onto approximation spaces}

Let us first introduce three $L^2$-projections on the spaces we have
defined. First, consider the projection $\ProjpK: L^2(K)\to \VpK$ on the
polynomial modes: for any $f\in L^2(K)$, find  $\pi^K_p f\in \VpK$ such that
\begin{align}
	( \ProjpK f, v_p)_K  = (f,v_p)_K, \qquad\forall v_p \in \VpK. \label{projp}
\end{align}
Next, define the projection $\ProjnK: L^2(K)\to \VnK$ on the
piecewise constant sub-grid space:  for any $f\in L^2(K)$, find  $\pi^K_n f\in \VnK$ such that
\begin{align}
	( \ProjnK f, v_n)_K  = (f,v_n)_K, \qquad\forall v_n \in \VnK. \label{projn}
\end{align}
Finally, define the projection on the combined space $\pi^K_\delta: L^2(K)\to V_\delta(K)$ given by: for any $f\in L^2(K)$, find  $\pi^K_\delta f\in V_\delta(K)$ such that
\begin{align}
	( \pi^K_\delta f, v_\delta)_K  = (f,v_\delta)_K, \qquad\forall v_\delta \in V_\delta(K). \label{projK}
\end{align}
Note that since the spaces $\VpK$ and $\VnK$ are in general not orthogonal with respect to $L^2(K)$, the projection $\pi_\delta^K$ is not simply the sum of $\ProjpK$ and $\ProjnK$. Also, all of these projections can be defined on the global spaces in a straight-forward way. The elementwise extension of $\pi_\delta^K$ to $\Omega$ is simply denoted by $\pi_\delta$.
Figure~\ref{fig:projection} illustrate the projection of $\pi^K_\delta f$ of two functions $f$ as well as the components of $\pi^K_\delta f$ lying in $\VnK$ and $\VpK$.

Since the $L^2$-approximation is the best-approximation in the discrete space with respect to the $L^2$-norm and since the space $V_\delta(K)$ contains the polynomial space $\mathbb P_{\! p}(K)$ we obtain immediately the following estimate for the $L^2$-projection: Let $f\in H^k(K)$, then 
\[
	\| f-\pi^K_\delta f\|_K \le C \left( \frac{h_K}{p} \right)^s | f |_{H^s(K)}
\]
for all $1\le s \le \min(p+1,k)$ for some constant $C$ independent of $h_K, p$ and
$n$.  Similar results can be obtained for other Sobolev norms of the error,
see \cite{BookSchwab,BookQuarteroni}. 
Therefore, exponential convergence can be
obtained for analytic solutions under $p$-refinement.  On the other hand,
since the space contains the piece-wise functions on the sub-grid $T_K$, one
obtains the following estimate
\[
	\| f-\pi^K_\delta f\|_K \le C  \frac{h_K}{n} | f |_{H^1(K)}.
\]

Finally, we note that $\pi^K_\delta$ (and $\pi^K_n$) preserves local averages on the sub-cells $k_i^\ell$ since one can test in \eqref{projK} with each sub-cell piece-wise constant basis function.

\begin{figure}[t]
  \begin{center}
    \begin{minipage}{.24\textwidth}
      \begin{center}
        \includegraphics[width=\textwidth]{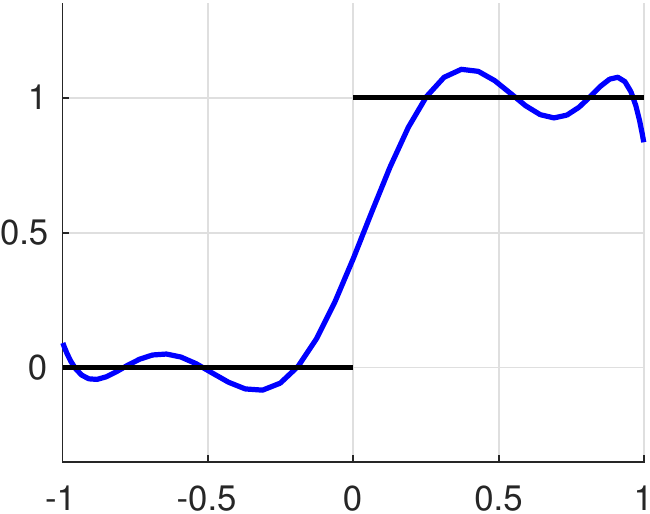} \\
        $n = 1$
      \end{center}
    \end{minipage}
    \begin{minipage}{.24\textwidth}
    \begin{center}
      \includegraphics[width=\textwidth]{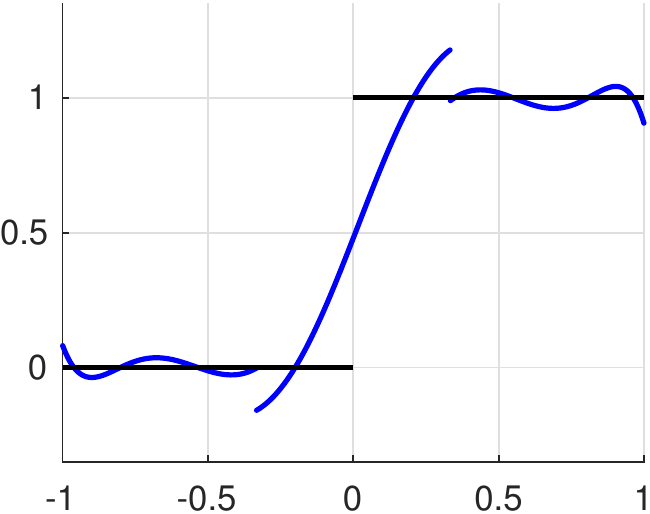} \\
      $n = 3$
    \end{center}
    \end{minipage}
    \begin{minipage}{.24\textwidth}
    \begin{center}
      \includegraphics[width=\textwidth]{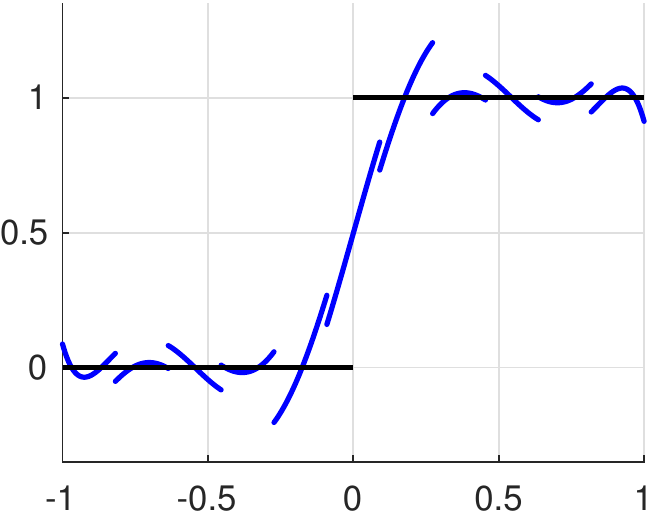} \\
      $n = 11$
    \end{center}
    \end{minipage}
    \begin{minipage}{.24\textwidth}
    \begin{center}
      \includegraphics[width=\textwidth]{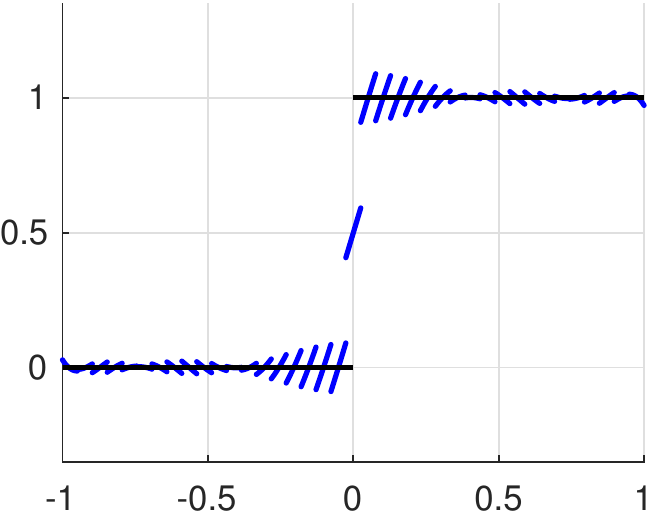} \\
      $n = 39$
    \end{center}
    \end{minipage}
    \\
    \begin{minipage}{.24\textwidth}
      \begin{center}
        \includegraphics[width=\textwidth]{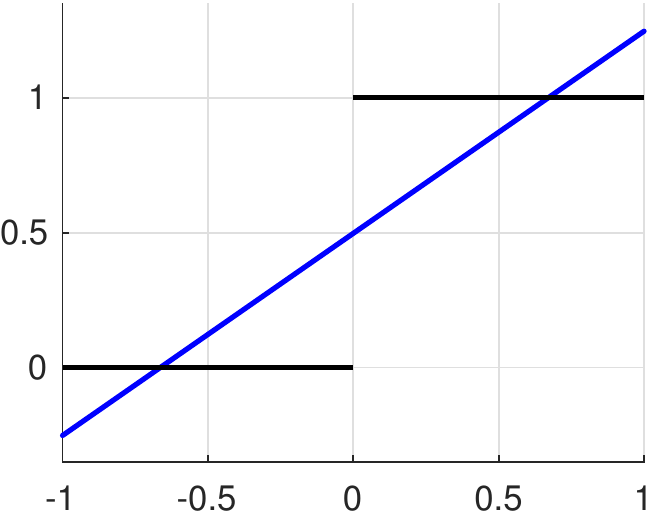} \\
        $p = 1$
      \end{center}
    \end{minipage}
    \begin{minipage}{.24\textwidth}
    \begin{center}
      \includegraphics[width=\textwidth]{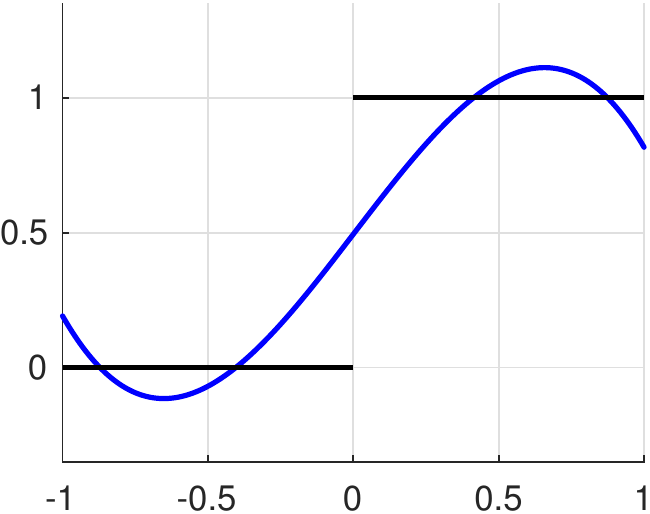} \\
      $p = 3$
    \end{center}
    \end{minipage}
    \begin{minipage}{.24\textwidth}
    \begin{center}
      \includegraphics[width=\textwidth]{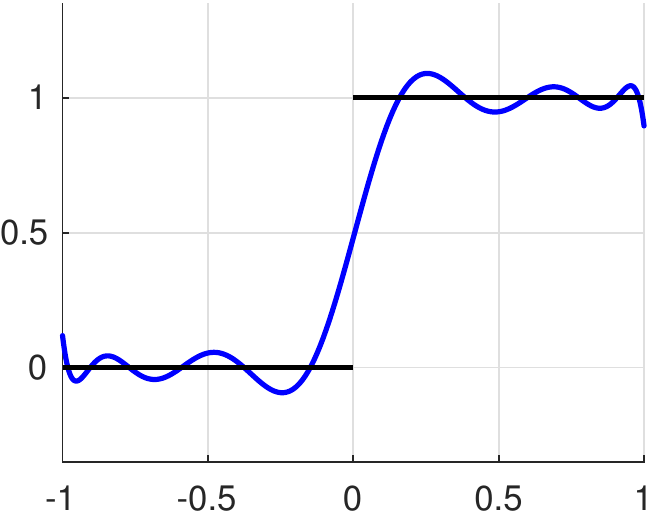} \\
      $p = 11$
    \end{center}
    \end{minipage}
    \begin{minipage}{.24\textwidth}
    \begin{center}
      \includegraphics[width=\textwidth]{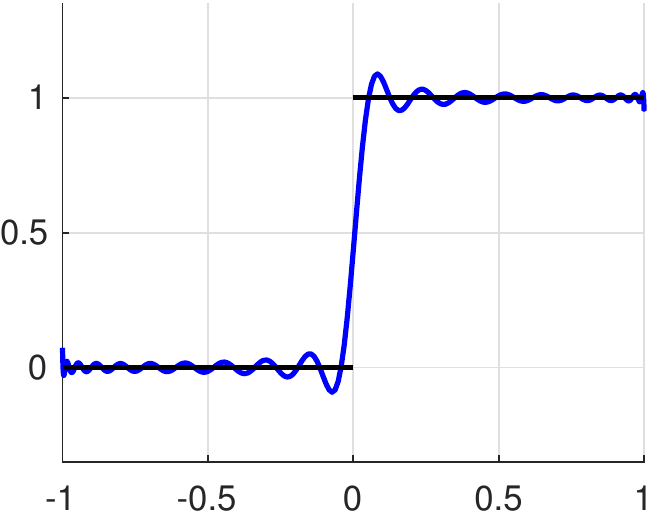} \\
      $p = 39$
    \end{center}
    \end{minipage} \\
  \end{center}
  \caption{$L^2$-approximation of a heaviside-function using the approximation space $V_\delta(K)$ introduced here with $\delta=(p=8,n)$ (top) and a polynomial approximation space (bottom), i.e. $\delta=(p,n=1)$. }
  \label{fig:overshoot}
\end{figure}

Figure \ref{fig:overshoot} illustrates the $L^2$-projection of the Heaviside-function onto the considered approximation space, as well as the projection onto a pure polynomial space.
As one can see, the presence of the piece-wise low order components does reduces the overshoot only for very (too) large $n$. 
However, in the asymptotic limit $n\to \infty$ we have observed that the overshoots are reduced with a rate of $1/n$ (the results are not reported here for the sake of compact presentation). 
Since this is not satisfactory, we explain in the next sub-section a strategy to project only on the low-order modes if necessary.

\subsection{Projections with penalized High-Order Components}

We note that the projections of monotonic functions on the $V_\delta$-space are in general not
monotonic, which can cause severe problems for the solution of non-linear
PDEs. Here, we propose a penalization approach in order to suppress the
polynomial modes within elements flagged by an indicator function.

Given a constant element penality parameter $\gamma\ge 0$ (see the upcoming Section~\ref{ssec:indicator} for an explicit construction), we extend the projection (\ref{projK}) by adding a penality term, to define the local penalized projection
$\pi_{\delta,\gamma}^K f$ of $f$ in element $K$ according to: for any $f\in L^2(K)$, 
find $\pi_{\delta,\gamma}^K f\in V_\delta(K)$ such that 
\begin{align}
  \big( \pi_{\delta,\gamma}^K f, v_\delta\big)_K +
  \gamma \, \big(\pi^K_p (\pi_{\delta,\gamma}^K f), \pi^K_p v_\delta\big)_K
  = \big(f,v_\delta\big)_K, \qquad\forall v_\delta \in V_\delta(K). \label{projKpenal}
\end{align}
This projection is well-defined since the bilinear form on the left hand side still defines a scalar product for any $
\gamma\ge 0$.
The motivation behind this projection is that the polynomial modes will
be continuously suppressed for increasing values of $\gamma$.
This can be easily seen from the equivalent minimization formulation for $\pi_{\delta,\gamma}^K f$:
\[
	\pi_{\delta,\gamma}^K f = \argmin_{w_\delta\in V_\delta(K)} 
	\Big[
	\big\| w_\delta -f\big\|_K^2 
	+
 	\gamma \, \big\|\pi^K_p (w_\delta)\big\|_K^2
  	\Big].
\]

\begin{figure}[t]
  \begin{center}
    \begin{minipage}{.24\textwidth}
      \begin{center}
        \includegraphics[width=\textwidth]{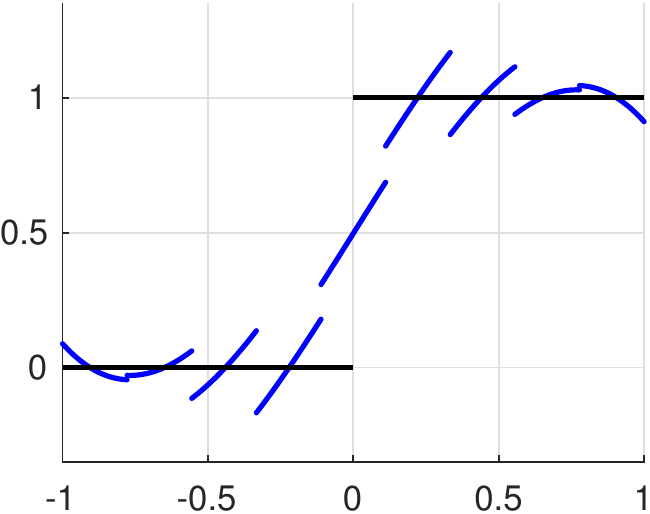} \\
      \end{center}
    \end{minipage}
    \begin{minipage}{.24\textwidth}
      \begin{center}
        \includegraphics[width=\textwidth]{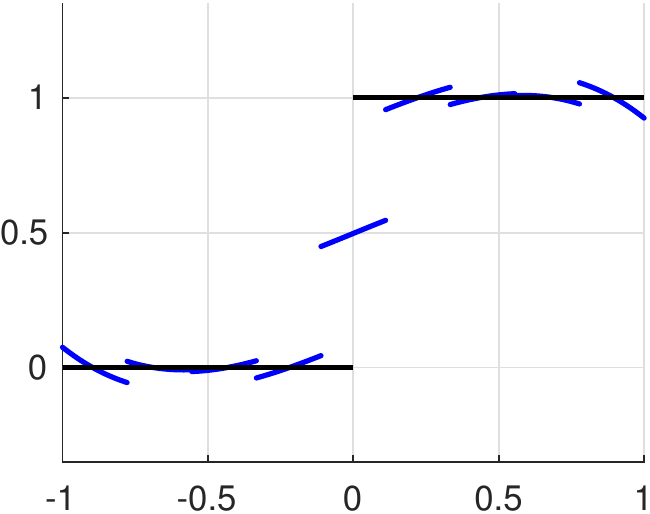} \\
    \end{center}
    \end{minipage}
    \begin{minipage}{.24\textwidth}
    \begin{center}
      \includegraphics[width=\textwidth]{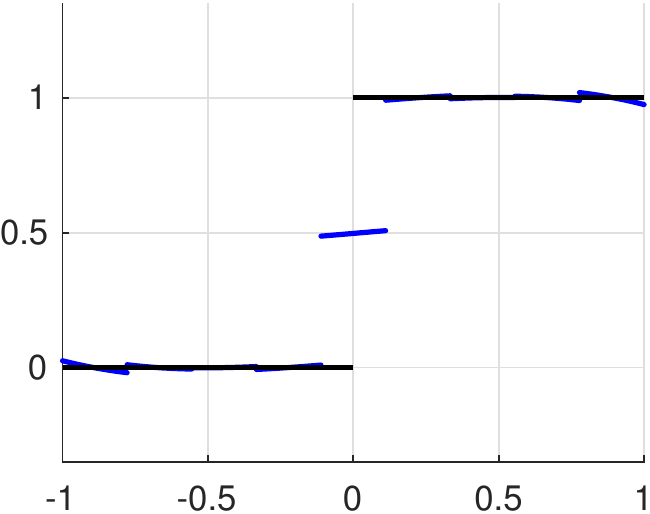} \\
    \end{center}
    \end{minipage}
    \begin{minipage}{.24\textwidth}
    \begin{center}
      \includegraphics[width=\textwidth]{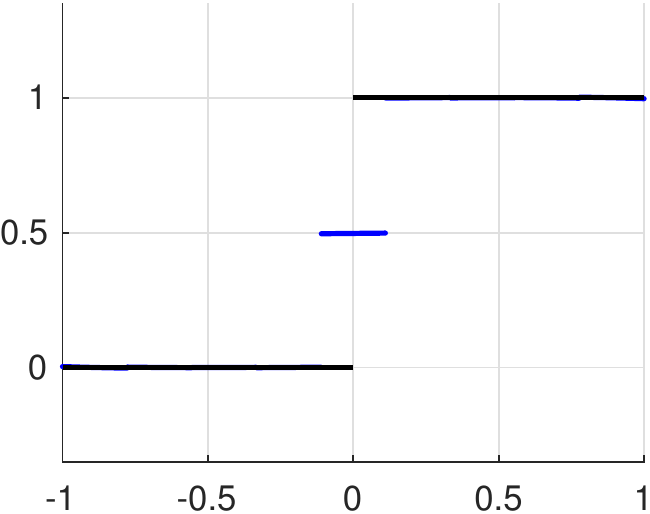} \\
    \end{center}
    \end{minipage}
    \\
    \begin{minipage}{.24\textwidth}
      \begin{center}
        \includegraphics[width=\textwidth]{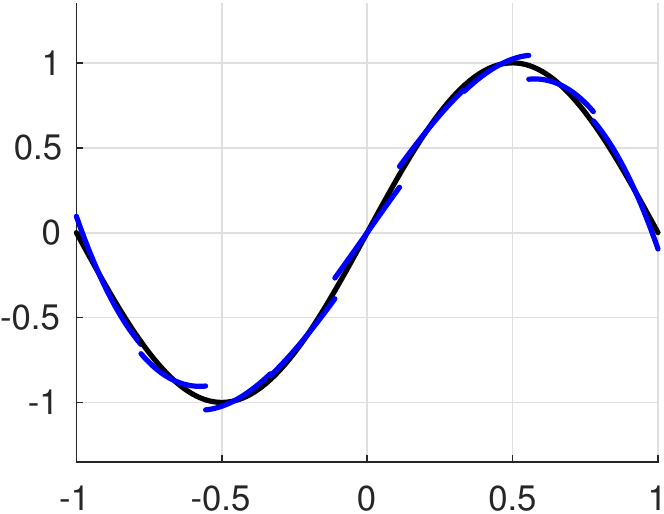} \\
        $\gamma = 0$
      \end{center}
    \end{minipage}
    \begin{minipage}{.24\textwidth}
    \begin{center}
      \includegraphics[width=\textwidth]{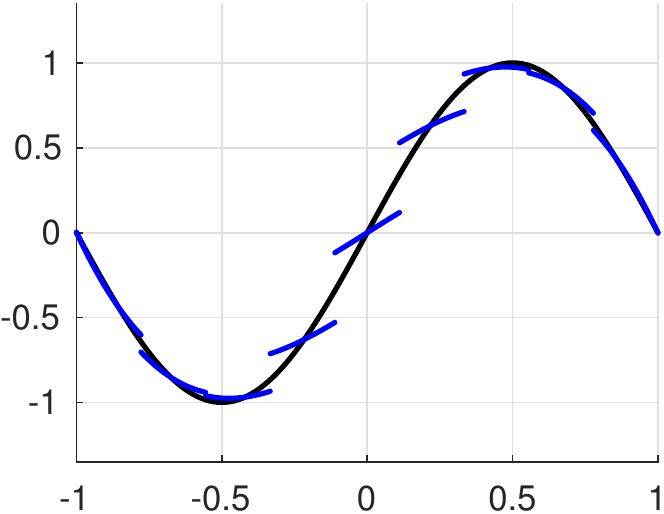} \\
      $\gamma = 0.05$
    \end{center}
    \end{minipage}
    \begin{minipage}{.24\textwidth}
    \begin{center}
      \includegraphics[width=\textwidth]{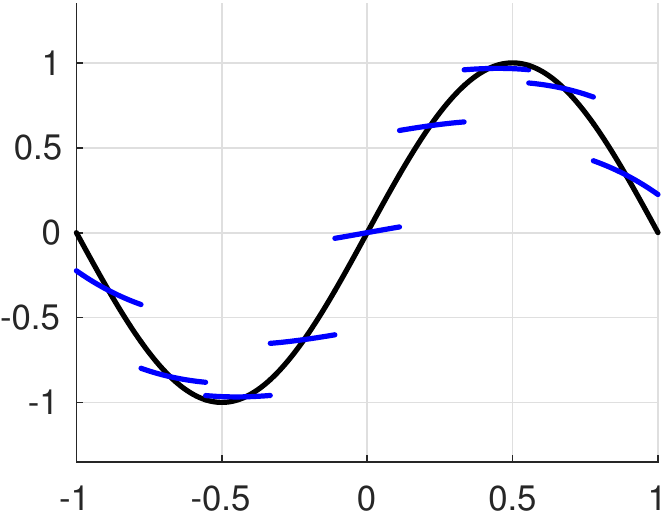} \\
      $\gamma = 0.5$
    \end{center}
    \end{minipage}
    \begin{minipage}{.24\textwidth}
    \begin{center}
      \includegraphics[width=\textwidth]{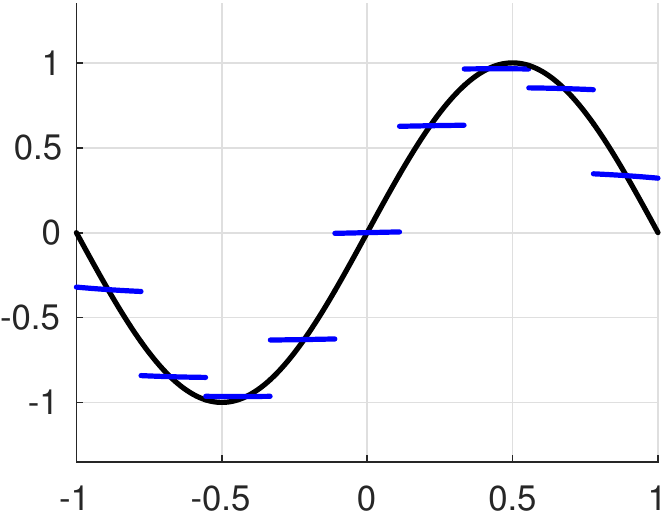} \\
      $\gamma = 5$
    \end{center}
    \end{minipage}
  \end{center}
  \caption{Penalty-based $L^2$-projection ($p=4$, $n=9$) for an increasing sequence of $\gamma$ of the Heaviside-function (top) and the sin-function (bottom).}
  \label{fig:penalproj}
\end{figure}

In Figure \ref{fig:penalproj} we present the penalty-based $L^2$-projection for a sequence of increasing values of $\gamma$ of the Heaviside-function as well as a smooth sin-function. As we can see, as is theoretically justified, the result converges to the monotonic $L^2$-projection using only the low-order modes.

\subsection{Indicator}
\label{ssec:indicator}
Having an indicator for the penalty-based $L^2$-projection at hand is an essential ingredient. 
A number of different sensors or indicators have been proposed in the
literature, to identify so-called ``trouble cells'' where some kind of
additional stabilization procedure is required to avoid oscillations
in the solution. In our setting, we find it natural to define such an
indicator directly in terms of the spaces and projectors we use.

In particular, we define a sensor function $s_K$ that measures {\it the ability to represent
a function $u_\delta$ within an element $K$ by a pure polynomial function
while still maintaining the local average on the sub-cells.} 
Now, the pure polynomial function that we mention above could in practice be taken as the $L^2$-projection \eqref{projp} onto $\widehat V_p(K) := \VpK\oplus \mbox{span}\{1\}$, i.e. where in addition $\VpK$ is augmented with the constant functions in order to obtain the full polynomial space.
However, this projection is not the best polynomial function that minimizes the local average of the original function. 
This motivates to introduce another projection $\overline{\pi}_{p,n}^K:V_\delta (K) \to \hVpK$ given by: 
for any $u_\delta\in V_\delta(K)$, 
find $\overline{\pi}_{p,n}^K u_\delta\in \hVpK$ such that 
\begin{align}
	\langle  \overline{\pi}_{p,n}^K u_\delta, v_p\rangle_{K}  
	= 
	\langle u_\delta ,v_p\rangle_{K}, 
	\qquad\forall v_p \in \hVpK, \label{projpn}
\end{align}
where $\langle \cdot, \cdot \rangle_{K}  $ is defined by 
\begin{align*}
	\langle w_\delta, v_\delta\rangle_{K}  
	:= 
	( \ProjnK w_\delta,\ProjnK v_\delta)_K,
	\qquad\forall w_\delta,v_\delta \in V_\delta(K).
\end{align*}
Indeed, the projection $\overline{\pi}_{p,n}^K u_\delta$ is the solution to the following minimization problem
\[
	\inf_{v_p\in \hVpK} \| \ProjnK (u_\delta - v_p) \|_K.
\]
Of course, this projection is in general not well-posed, in particular if $N(p)+1 > n$, i.e. if there are more degrees of freedom than conditions so that uniqueness is not given. 
\begin{lemma}
	\label{lem:solv}
	Consider a uniform partition of the unit simplex into $n=(r+1)^d$ simplices of edge-length $1/(r+1)$. If $r\ge p$, then the mapping $\ProjnK: \hVpK \to  \VnK$ is injective. 
\end{lemma} 

The proof is presented in the appendix.

\begin{remark}
	Numerical results indicate that Lemma~\ref{lem:solv} does not provide a sharp condition on the solvability for large values of $p$. For example, for $d=2$ and $p=4$, the mapping is injective even for $r=3$ which generates a sub-grid of 16 triangles whereas $N(4)+1 = 15$.
	In practise, it is easy to check solvability on the unit simplex by numerical computation for any given parameters $p$ and $n$ in given dimension $d$ once and for all.
\end{remark}

It is easy to see that the injectivity of $\ProjnK:\hVpK \to \VnK$ implies that $\| \ProjnK \cdot \|_K$ is a norm on $\VpK$ and that the  projection $\overline{\pi}_{p,n}^K$ defined by \eqref{projpn} is in consequence uniquely determined.

Now, having clarified the well-posedness of the projection $\overline{\pi}_{p,n}^K$, we are now ready to state the indicator function
\begin{align}
	 s_K =  \| \ProjnK ( u_\delta -\overline{\pi}_{p,n}^K u_\delta ) \|_{\infty,K},
\end{align}
with $\| \cdot \|_{\infty,K}$ being the $L^\infty(K)$-norm.
Therefore, the sensor $s_K$ can be viewed as the maximal error if one would replace the current approximation $u_\delta$ by the  (pure) polynomial that best preserves the local averages on the sub-cells.

We also define a local normalization factor for this sensor, which is given by
\begin{align}
 s^0_K =  \| \ProjnK  u_\delta  \|_{\infty,K} + s_\varepsilon,
\end{align}
where $s_\varepsilon$ is a small number to avoid zero division.


In terms of this sensor $s_K$, we now define the penalty function
\begin{align}
  \gamma_K = C_\mathrm{pen} \max \left( 0, \frac{s_K}{s^0_K} - \tau \right)
\end{align}
where $C_\mathrm{pen}$ is a problem and discretization dependent parameter, which we
empirically set to $C_\mathrm{pen} = 10^{7}$ in all our examples. 
The threshold $\tau$ is used to ensure zero penalization for smooth solutions and small perturbations, and we
set it empirically to $\tau =0.01 / p$.

\begin{figure}[t]
  \begin{center}
    \begin{sideways} \hspace{-8mm}Function $u_\delta$ \end{sideways}
    \begin{minipage}{.23\textwidth}
      \begin{center}
        \includegraphics[width=\textwidth]{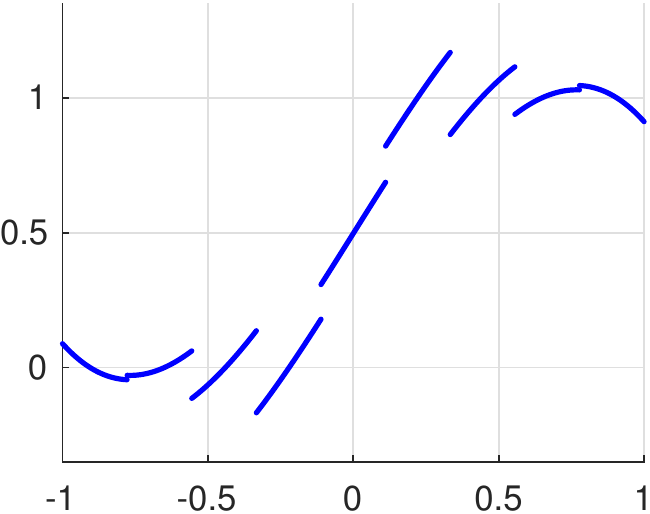} \\
      \end{center}
    \end{minipage}
    \begin{minipage}{.23\textwidth}
      \begin{center}
        \includegraphics[width=\textwidth]{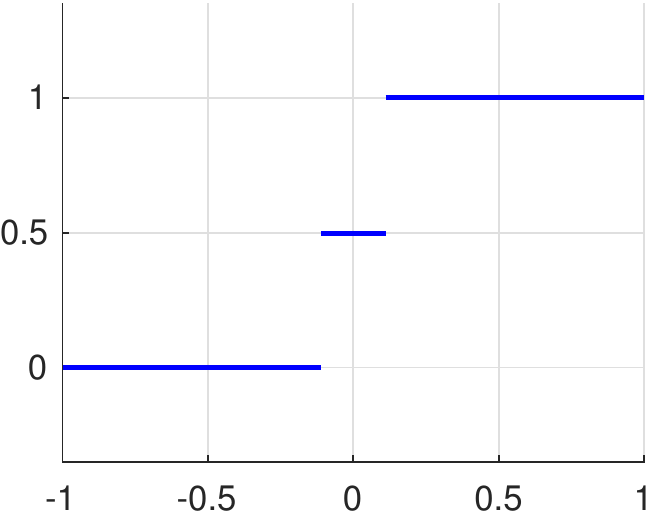} \\
    \end{center}
    \end{minipage}
    \begin{minipage}{.23\textwidth}
    \begin{center}
        \includegraphics[width=\textwidth]{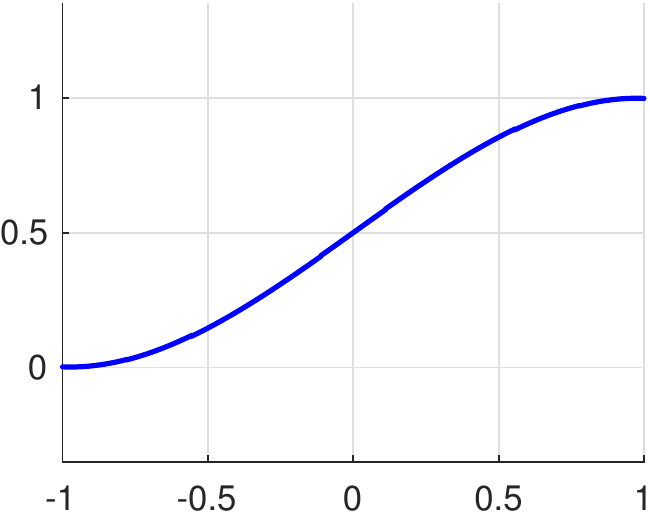} \\
    \end{center}
    \end{minipage}
    \begin{minipage}{.23\textwidth}
    \begin{center}
        \includegraphics[width=\textwidth]{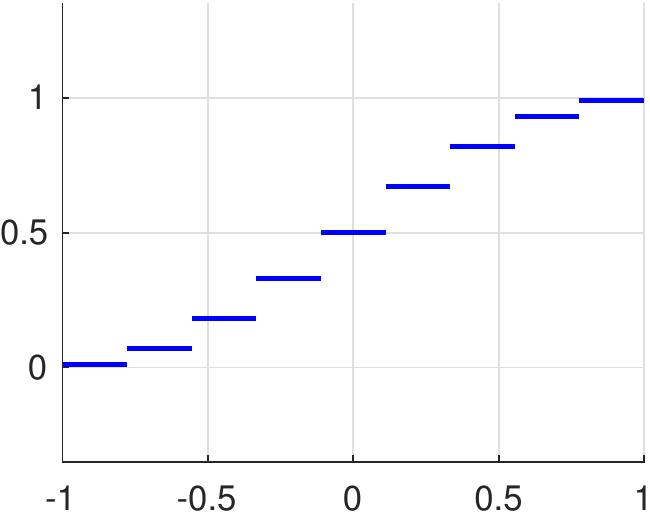} \\
    \end{center}
    \end{minipage}
    \\
    \begin{sideways} \hspace{-6mm}Sensor  \end{sideways}
    \begin{minipage}{.23\textwidth}
      \begin{center}
        \includegraphics[width=\textwidth]{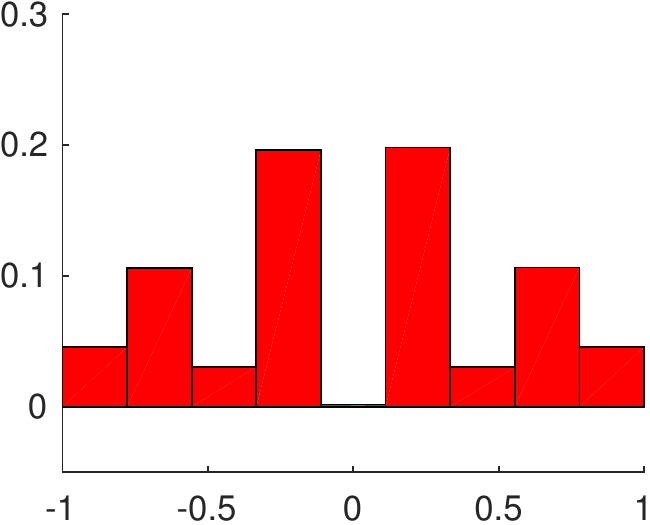} \\
      \end{center}
    \end{minipage}
    \begin{minipage}{.23\textwidth}
    \begin{center}
        \includegraphics[width=\textwidth]{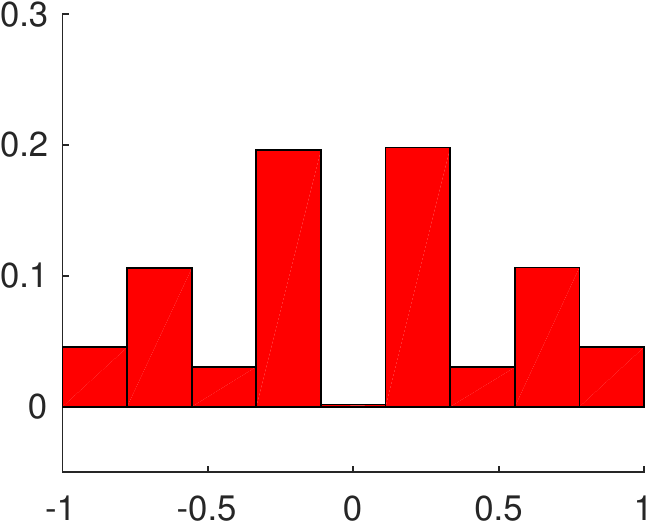} \\
    \end{center}
    \end{minipage}
    \begin{minipage}{.23\textwidth}
    \begin{center}
        \includegraphics[width=\textwidth]{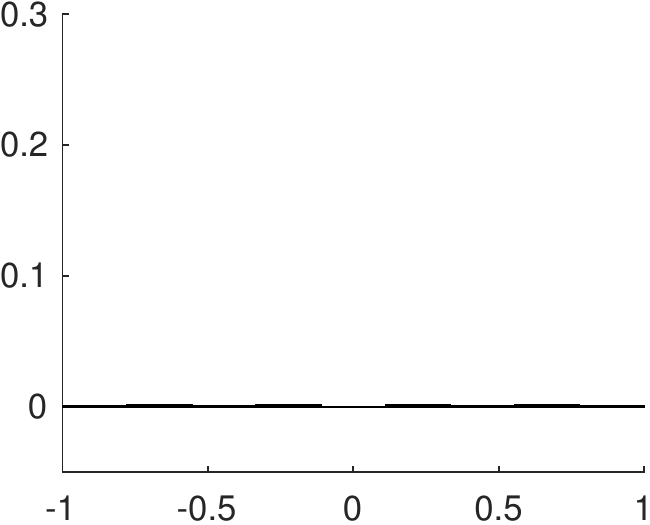} \\
    \end{center}
    \end{minipage}
    \begin{minipage}{.23\textwidth}
    \begin{center}
        \includegraphics[width=\textwidth]{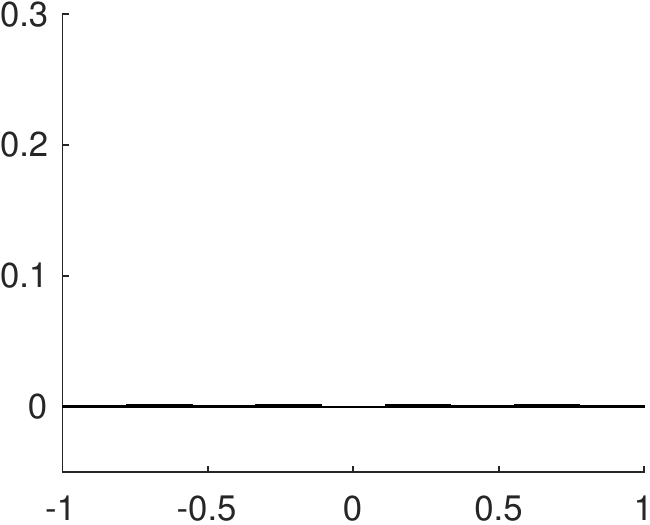} \\
    \end{center}
    \end{minipage}
  \end{center}
  \caption{Illustration of the sensor function $|\ProjnK ( u_\delta -\overline{\pi}_{p,n}^K u_\delta )|$ over one element $K=(-1,1)$ (bottom) for $p=4$, $n=9$ applied to four different functions (top).} 
  \label{fig:sensorexample}
\end{figure}

Figure \ref{fig:sensorexample} illustrates the sensor function $|\ProjnK ( u_\delta -\overline{\pi}_{p,n}^K u_\delta )|$ for four different functions using the discretization parameters $p=4$ and $n=9$. 
We observe that the sensor indicates the first two functions as they can not be represented accurately by the $L^2$-projection onto polynomials only. 
We finalize this section with listing some properties of the sensor:
\begin{itemize}
\setlength{\itemsep}{0pt}
\item If $u_\delta\in \hVpK$, then $s_K = 0$ since $\overline{\pi}_{p,n}^K u_\delta = u_\delta$. Figure \ref{fig:sensorexample} (third column) illustrates this property: this function is a pure polynomial function (with $p=4$) so that $\overline{\pi}_{p,n}^K u_\delta=u_\delta$. 
The relatively small error should not impact the sensor which motivates a positive threshold parameter $\tau = 0.01/p$.
\item If $u_\delta\in \VnK$, then $s_K$ is in general positive. It can nevertheless be zero or very small if the function can be nicely represented by a polynomial while keeping the local averages in the sub-cells, which is an important feature to allow for recovering the high-order representation, e.g. after a shock-transition. This is illustrated in Figure  \ref{fig:sensorexample} (fourth column) which does not sense this function (although it is a piecewise constant function).
\end{itemize}

\section{Discontinuous Galerkin formulation for conservation laws}

We devote now our attention to the discretization of conservation laws using the previously introduced approximation space~$V_\delta$.

\subsection{Governing equations}
\label{ssec:GovEqns}

Consider a general system of first-order conservation laws
\begin{align}
  \frac{\partial u}{\partial t} + \nabla \cdot F(u) = 0 \qquad
  \text{ in } \Omega,
\end{align}
with prescribed initial condition $u(x,t=0)=u_0$ and appropriate
boundary conditions on $\partial \Omega$ (more details will be provided in the results
section). 
Here $u=(u_1,\ldots,u_m)$ is the solution vector function and we will
consider the following equation:
\begin{description}
[leftmargin=.75cm,itemsep=10pt,font=\mdseries]
\item[\bf {\tt [1]} Convection:] 
	We consider $m=1$ and $F(u) = \beta u$ for a given velocity field
  $\beta$ (possibly space and time varying).
\item[\bf {\tt [2]} Inviscid Burgers' equation:] 
	The spatial dimension is limited to $d=1$ only; with $m=1$ and  $F(u)=u^2/2$.
\item[\bf {\tt [3]} Euler's equations of gas dynamics:] 
Here, $m=d+2$ and the system of conservation laws is given by
\begin{align}
\frac{\partial \rho}{\partial t}  + \frac{\partial}{\partial x_i}
(\rho u_i) &= 0, \label{ns1} \\
\frac{\partial}{\partial t} (\rho u_i) +
\frac{\partial}{\partial x_i} (\rho u_i u_j+ p)  &= 0
\quad\text{for }i=1,2,3, \label{ns2} \\
\frac{\partial}{\partial t} (\rho E) +
\frac{\partial}{\partial x_i} \left(u_j(\rho E+p)\right) &= 0 \label{ns3}
\end{align}
where $u_i$, $i=1,\ldots,d$ are the velocity components, $\rho$ is the fluid density and $E$ is the total energy, thus the unknown vector function consists of $u=(u_1,u_2,u_3,\rho,E)$.
Further, we assume an ideal gas with
pressure $p$ of the form
\begin{align}
p=(\gamma_{\mathsf{a}}-1)\rho \left( E - \frac12 u_k u_k\right),
\end{align}
where $\gamma_{\mathsf{a}}$ is the adiabatic gas constant.
  
\end{description}

\subsection{Spatial Discretization}

Our discretization is a standard (discontinuous) Galerkin form on our space
$V_\delta$, with the optional penalization of the polynomial modes as
described above. Define the multi-component spaces $V^m_{\delta}(K) = [
  V_\delta(K) ]^m$ and $V^m_{\delta} = [ V_\delta ]^m$. 
For all $t\in(0,T)$,
find $u_\delta(t)=(u_{\delta,1}(t),\ldots,u_{\delta,m}(t)) \in V^m_\delta$
such that
\begin{align}
  \int_K \partial_t u_\delta(t)\,v_\delta\
  - \sum_{k\in T_K} \int_k F\big(u_\delta(t)\big) \cdot \nabla v_\delta
  + \sum_{k\in T_K} \oint_{\partial k} \widehat{F}\big(u^+_\delta(t),u^-_\delta(t),n_k\big)\,v_\delta \\ \nonumber
  + \gamma_K{(u_\delta(t))} \int_K \left(\ProjpK u_\delta(t)\right)(\ProjpK v_\delta)   = 0,
  \qquad \forall v_\delta \in V^m_\delta(K), \label{femform}
\end{align}
on each element $K\in T$. We are using a compact notation to handle the
multi-component system.
Here, $\widehat{F}\left(u^+_\delta(t),u^-_\delta(t),n_k\right)$ is a numerical flux
function involving the solution on each side of the boundary $\partial k$ and
the outward normal vector $n_k$.

Next, write the global $u_\delta$ approximation in terms of the basis functions $\phi_i$:
\begin{align}
u_\delta(t) = \sum_{\ell=1}^{N_T} \sum_{i=1}^{N(p)+n} U^\ell_i(t) \, \varphi^\ell_i
\end{align}
where $U$ is the vector of all components $U^\ell_i=(U^\ell_{i,1},\ldots,U^\ell_{i,m})$, $\ell=1,\ldots,N_T$,
$i=1,\ldots,N(p)+n$. 
Note that the dimension of $U$ is $m\times N_T\times (N(p)+n)$.
Imposing equation~\eqref{femform} for each basis function leads to a semi-discrete system of the form
\begin{align}
\rM \, \dot{U}(t) + \Gamma(U(t)) \, \rM_{pp} \, U(t) = \rR(U(t)) \label{ODE}
\end{align}
Here, $\rM$ is a block-diagonal mass matrix, $\rM_{pp}$ is a (singular) block-diagonal
mass matrix for the polynomial components only, $\Gamma(U(t))$ is a diagonal
matrix with the values of the penalty parameter $\gamma_K$ for each element,
and $\rR(U(t))$ contains the remaining terms (which is in general non-linear).

The initial condition, denoted by $U_0$, is obtained by projecting the initial function $u_0$ onto the discretization space~$V_\delta$.

\subsection{Time integration by Implicit-Explicit Runge-Kutta methods}

The system (\ref{ODE}) can be integrated by a number of different ODE
solvers, including fully implicit and fully explicit or combined implicit-explicit (IMEX) schemes. 
Since the penalty term $\Gamma(U(t)) \, \rM_{pp} \, U(t)$ is very stiff we apply an IMEX scheme using the splitting
\[
	\rM \, \dot{U}(t) = f(U(t)) + g(U(t)),
\]
with explicit part $f$ and implicit part $g$ given respectively by
\[
	f(U(t)) = \rR(U(t)),
	\qquad\qquad
	g(U(t)) = -\Gamma(U(t)) \, \rM_{pp} \, U(t).
\]
In particular, we use the approach presented in \cite{persson13shock} where we freeze however the penalty function $\Gamma$ during the different stages. One time-step from $U_n\approx U(t^n)$ to $U_{n+1} \approx U(t^{n+1})$ reads as:
\vspace{5pt}
\begin{center}
 \fcolorbox{black}{shadecolor}{
\begin{minipage}{0.98\textwidth}
\begin{flushleft}
Set $\Gamma_n = \Gamma(U_{n})$ \\
\textbf{for} 			$i=1$ \textbf{to} $s$ \\
\hspace{8mm} 		$\displaystyle{U_{n,i} = U_{n} + \Delta t \sum_{j=1}^{i-1} \left [a_{i,j}\, r_j  +  \hat{a}_{i,j}\, \hat{r}_j\right]}$\\
\hspace{8mm} 		Evaluate $(\rM+\Delta t\, a_{i,i} \,\Gamma_n \, \rM_{pp}) \,{r}_{i}= - \Gamma_n \, \rM_{pp}  \, U_{n,i} $ \\
\hspace{8mm} 		Evaluate $\rM \,\hat{r}_{i}= \rR(U_{n,i} + \Delta t \,a_{i,i}\,r_{i})$ \\
\textbf{end for} \\
$\displaystyle{U_{n+1} = U_{n} + \Delta t  \sum_{i=1}^s\left[ b_j r_j +\hat{b}_j\; \hat{r}_j\right]}$.
\end{flushleft}
\end{minipage}
}
\end{center}
\vspace{5pt}
We use the ARS(2,2,2) scheme, with coefficients $a_{i,j}$, $\hat a_{i,j}$, $b_j$, $\hat b_j$ given by the matrices $A$, $\hat A$ and vector $b$, $\hat b$:
\[
	A =   
	\begin{pmatrix}
	0 & 0 & 0 \\
    0 & \alpha & 0 \\
	0 & 1-\alpha & \alpha 
  	\end{pmatrix}		
	\qquad
	\hat A =   
	\begin{pmatrix}
	0 & 0 & 0 \\
    \alpha & 0 & 0 \\
	\delta & 1-\delta & 0 
 	\end{pmatrix}		
 	\qquad
 	b =
	\hat b =   
	\begin{pmatrix}
	0 \\
    1-\alpha \\
	\alpha 
 	\end{pmatrix},		
\]
with fixed parameters $\alpha=1 - \frac{1}{\sqrt{2}}$, $\delta=-2\frac{\sqrt{2}}{3}$.

\section{Results}
We present here a collection of numerical results in one and two spatial dimensions with increasing complexity.
We use Roe's approximate Riemann solver for the numerical flux function in all our examples.

\subsection{Linear convection}
We start with the most simple case of scalar linear convection in one spatial dimension, see Section \ref{ssec:GovEqns} {\tt [1]}  with $\beta=1$, on $\Omega=(0,1)$ in order to illustrate some basic approximation properties. 
We consider two initial conditions $u_0=u(t=0)$, a Gaussian function and a Heaviside function, that are evolved under periodic boundary conditions for one cycle to obtain $u_\delta(t=1)$. 
The polynomial degree $p=1,2,3,4$ and the number of sub-cells $n=8$ per element is fixed while the number of elements $N_T$ is increased to generate convergence rates as $h\to 0$ where $h=1/n_T$ denotes the diameter of the elements.  
The time-step $\Delta t$ is kept small enough so that the error is dominated by the spatial discretization.

\begin{figure}
  \centering
  \begin{tabular}{cc}
    \includegraphics[width=.48\textwidth]{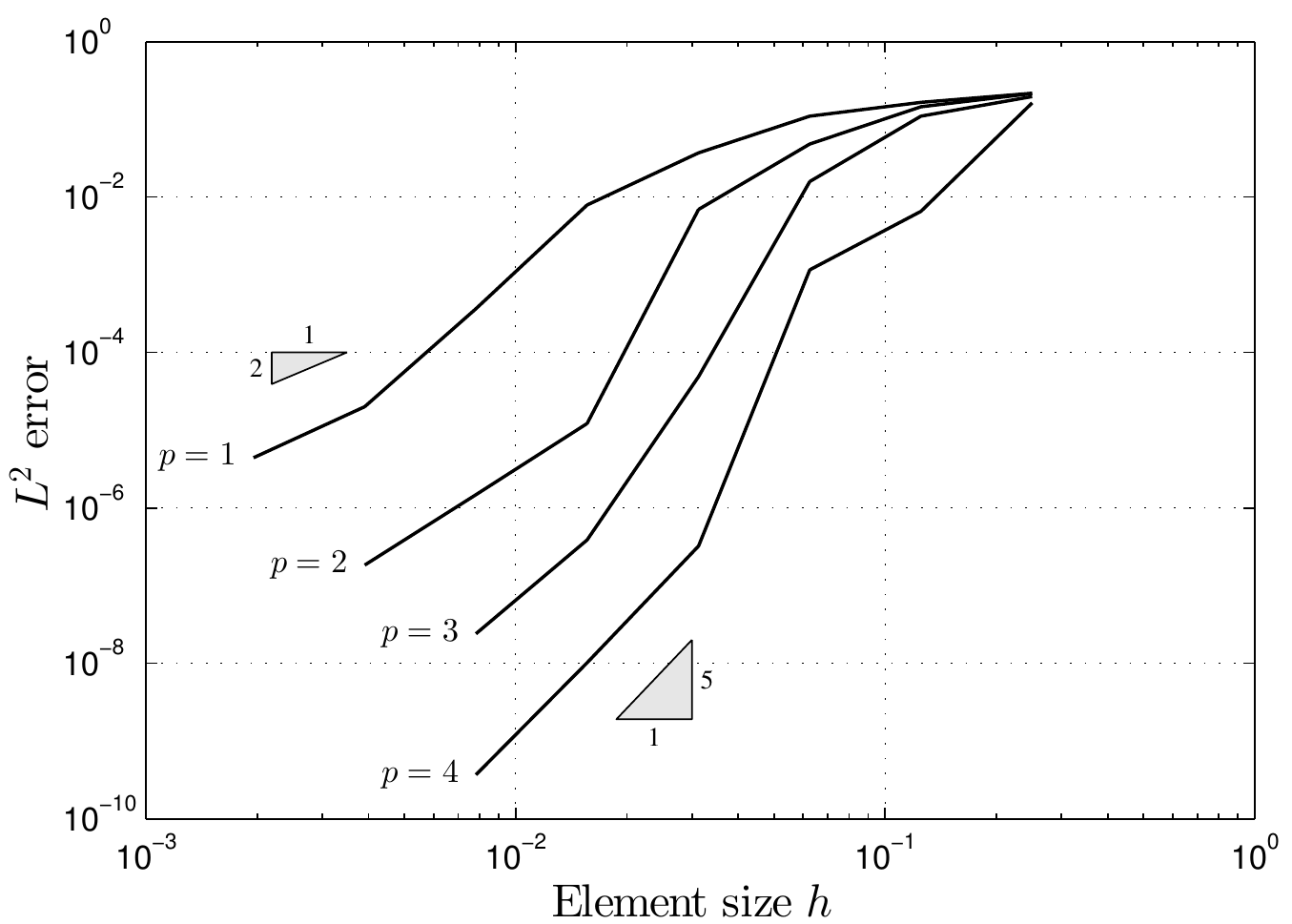} &
    \includegraphics[width=.48\textwidth]{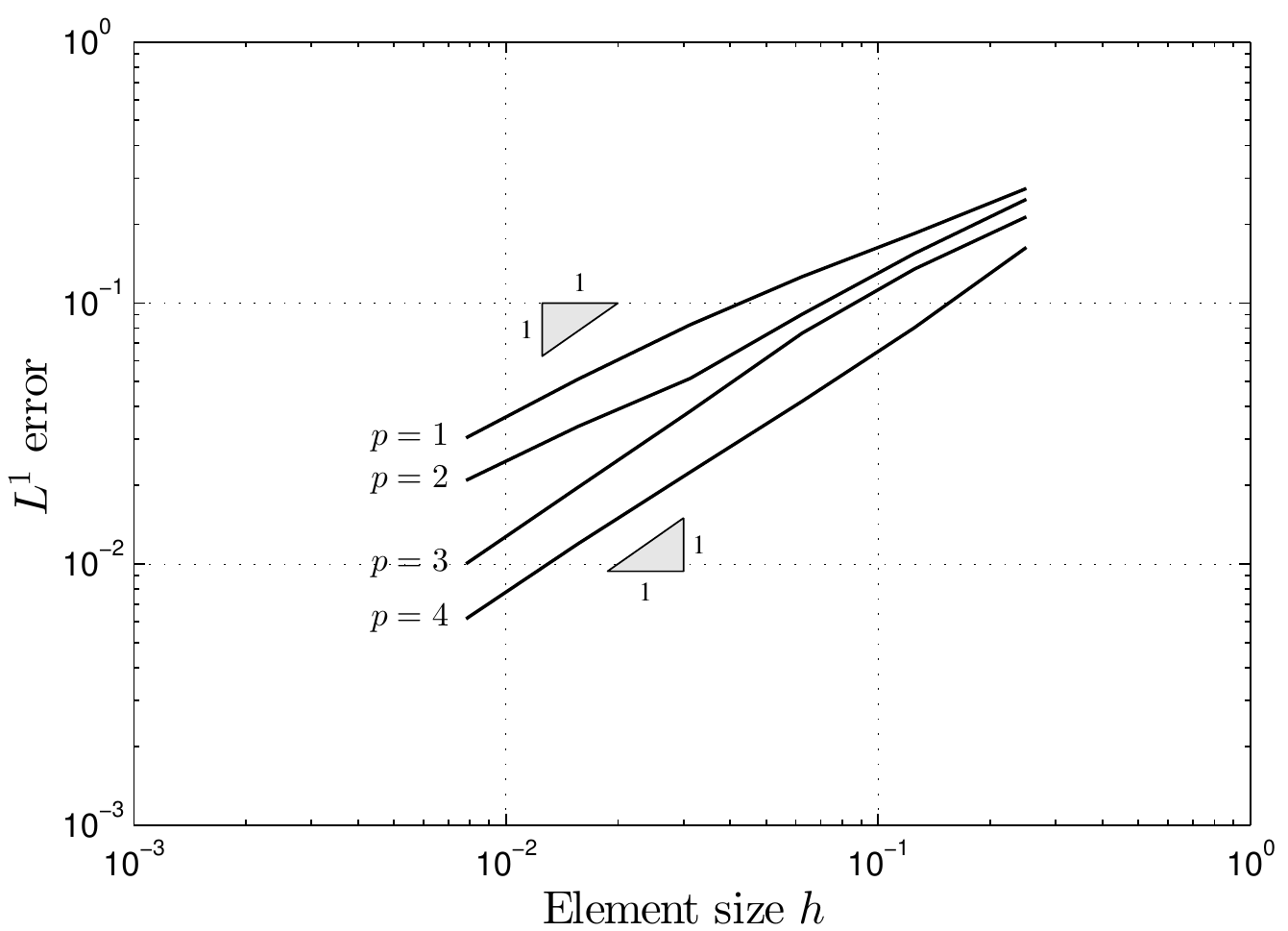} \\\\
    Gaussian initial condition & Heaviside  initial condition \\ 
  \end{tabular}
  \caption{Accuracy of the method for linear convection for a Gaussian (left) and Heaviside (right) initial condition with $p=1,2,3,4$ and  $n=8$ under mesh refinement at $t=1$.}
  \label{fig:convection}
\end{figure}

In Figure \ref{fig:convection} (left) we illustrate the $L^2(\Omega)$-norm of the error $\pi_\delta u_0 - u_\delta(t=1)$ for the different values of $p$ under mesh-refinement $h\to 0$ for the Gaussian function while \ref{fig:convection} (right) plots the same errors in the $L^1(\Omega)$-norm, which is standard for discontinuous solutions, for the Heaviside function.

\begin{figure}
  \centering
  \begin{tabular}{cccc}
    \includegraphics[width=.23\textwidth]{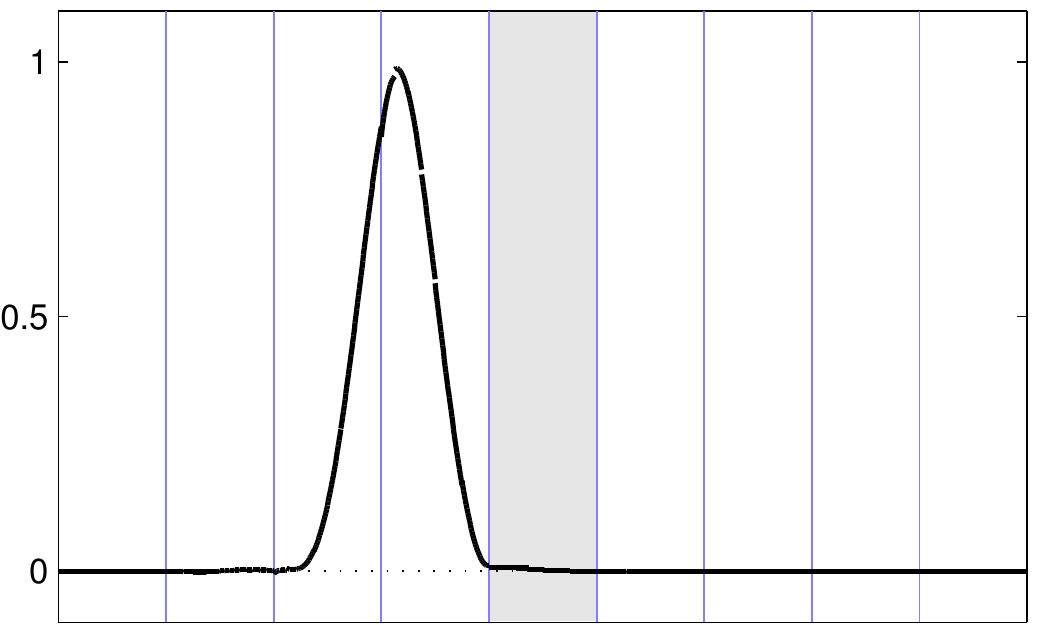} &
    \includegraphics[width=.23\textwidth]{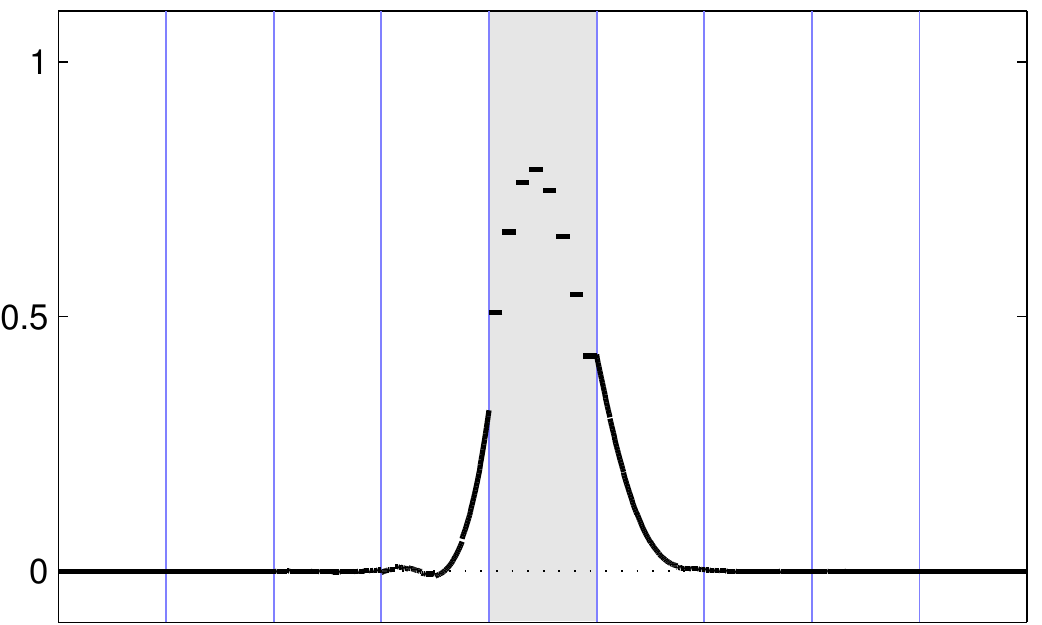} &
    \includegraphics[width=.23\textwidth]{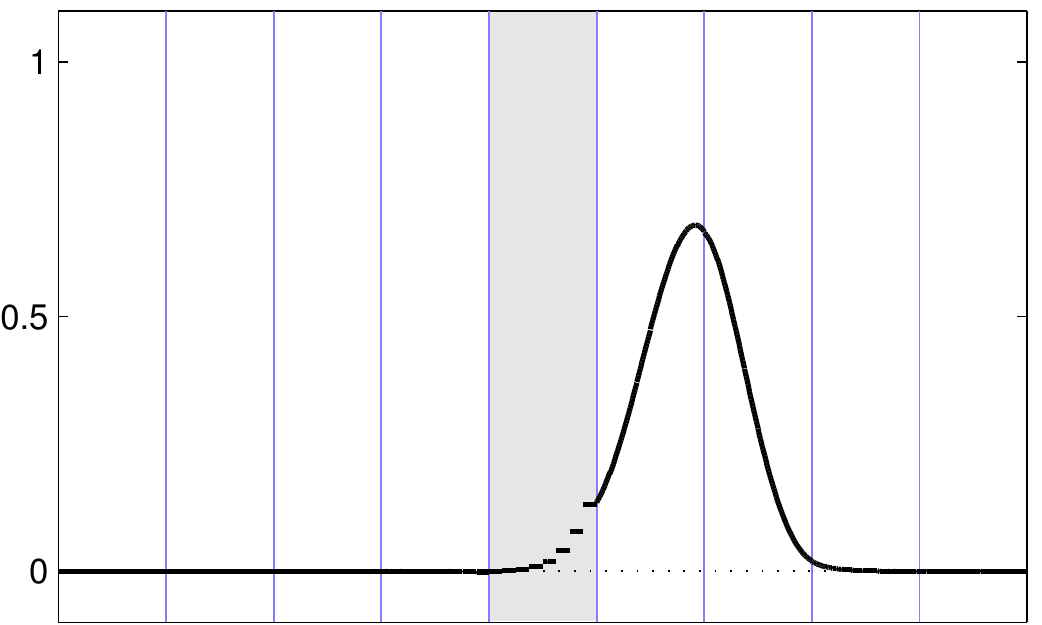} &
    \includegraphics[width=.23\textwidth]{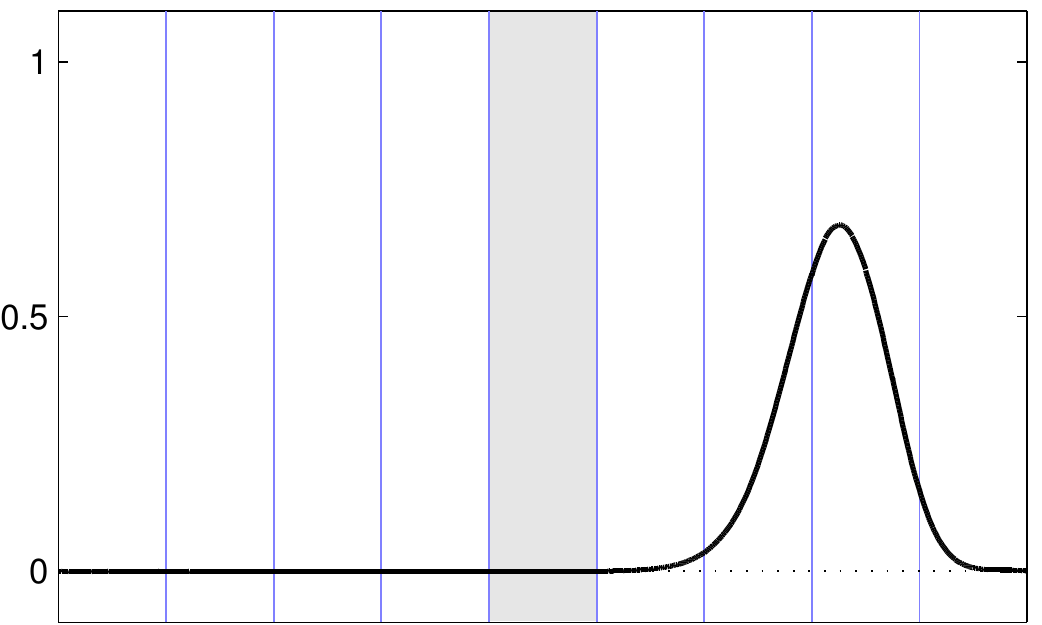} \\
  \end{tabular}
  \caption{Solution before, during and after it passes the marked element in which the penalization parameter is artificially activated.}
  \label{fig:recovery}
\end{figure}

Finally, in Figure~\ref{fig:recovery} we activate artificially the penalty parameter in one element and plot the solution of the transported Gaussian before, during and after it passes this element. We observe that the solution can qualitatively be recovered with the polynomial modes after the Gaussian traveled through the marked element, but that the shape suffers from the high dissipation that occurred during the piecewise constant representation in the marked element.

\subsection{Inviscid Burgers' equation}
\begin{figure}
  \centering
  \begin{tabular}{cccc}
    \includegraphics[width=.23\textwidth]{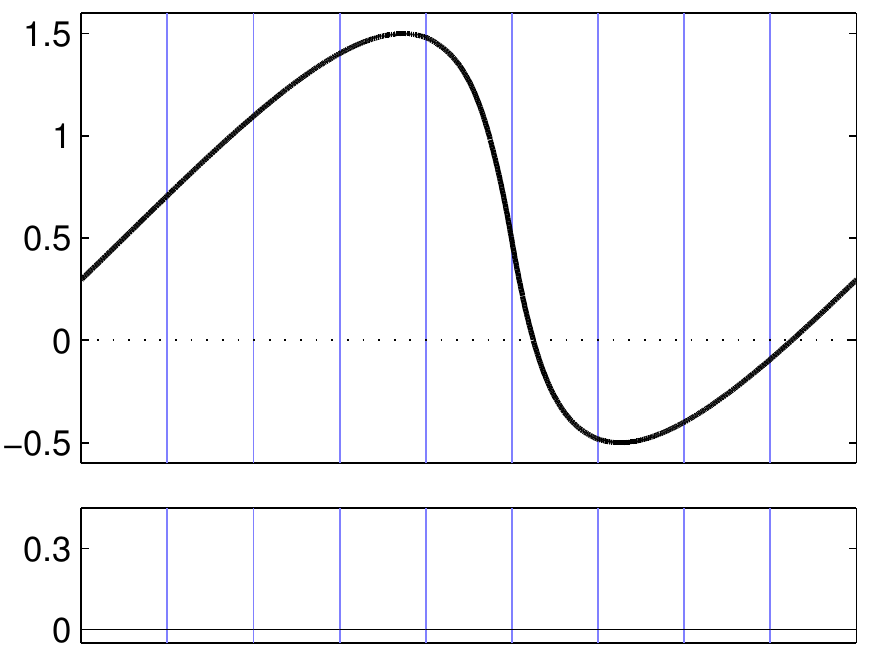} &
    \includegraphics[width=.23\textwidth]{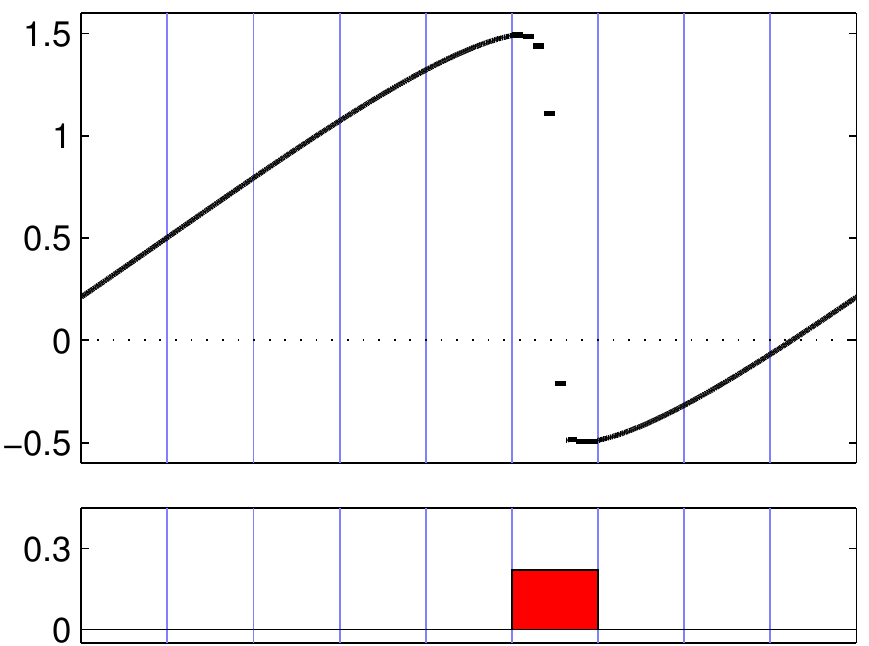} &
    \includegraphics[width=.23\textwidth]{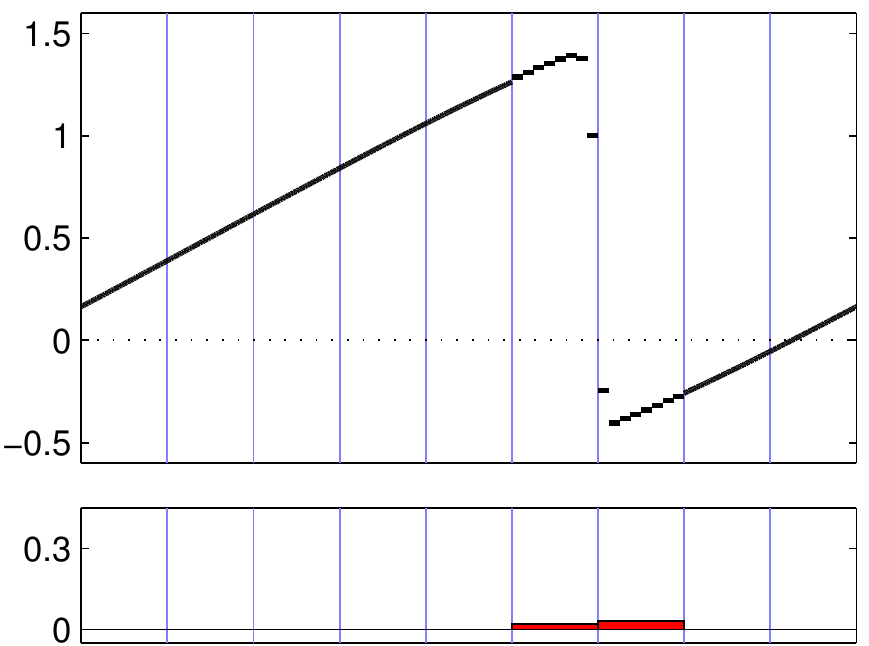} &
    \includegraphics[width=.23\textwidth]{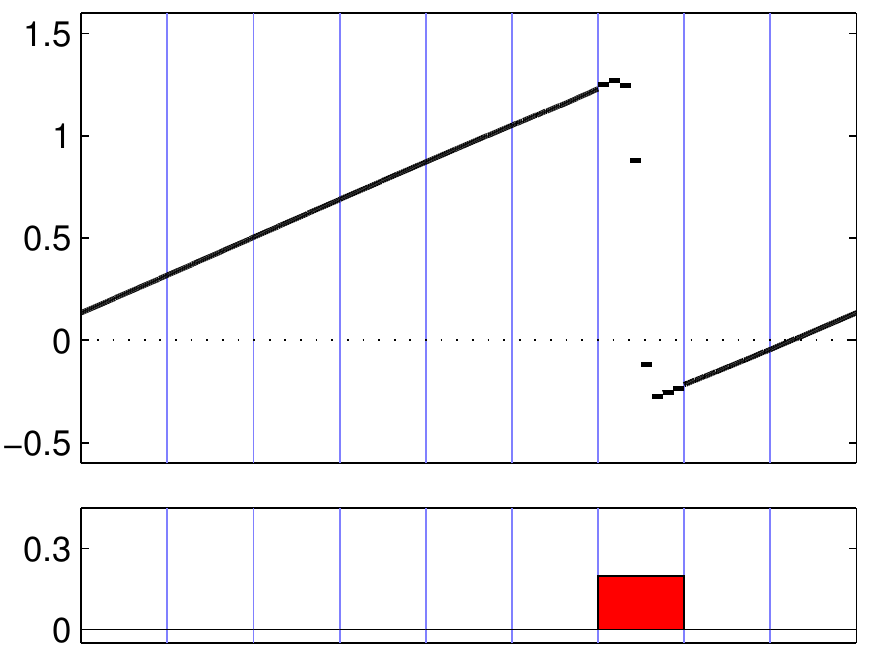} \\
    $t=0.11$ & $t=0.22$ & $t=0.33$ & $t=0.44$ \\ \\
    \includegraphics[width=.23\textwidth]{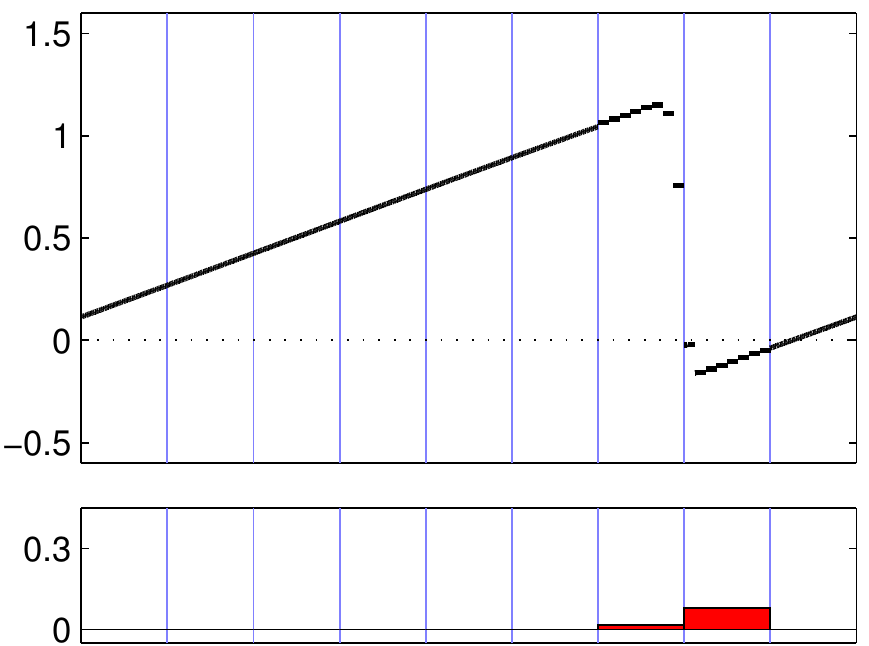} &
    \includegraphics[width=.23\textwidth]{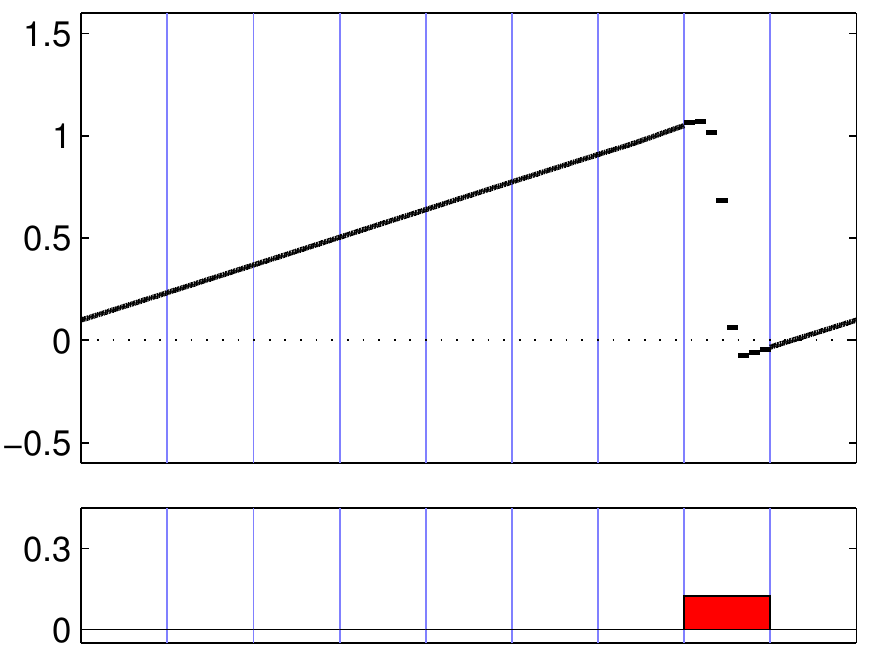} &
    \includegraphics[width=.23\textwidth]{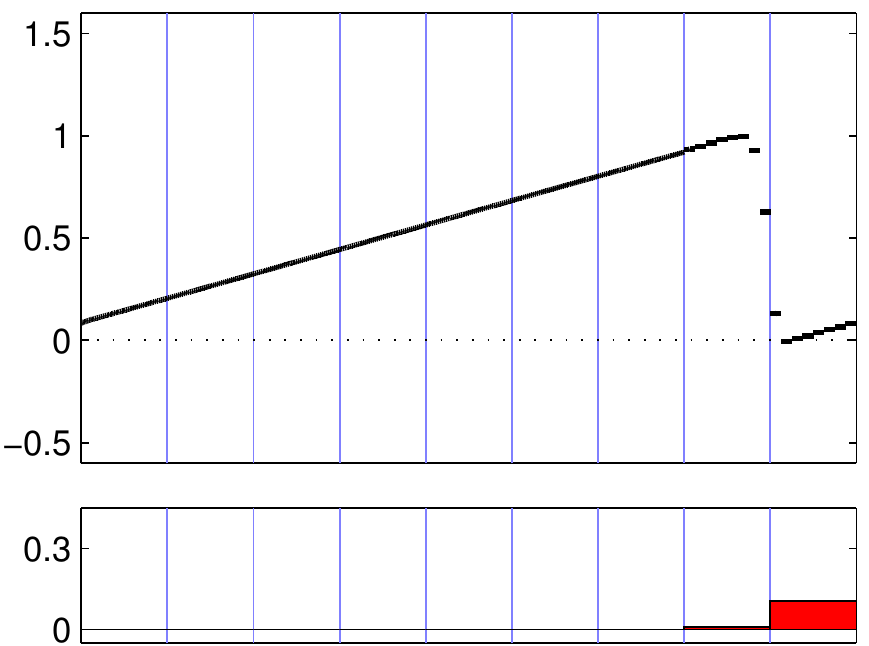} &
    \includegraphics[width=.23\textwidth]{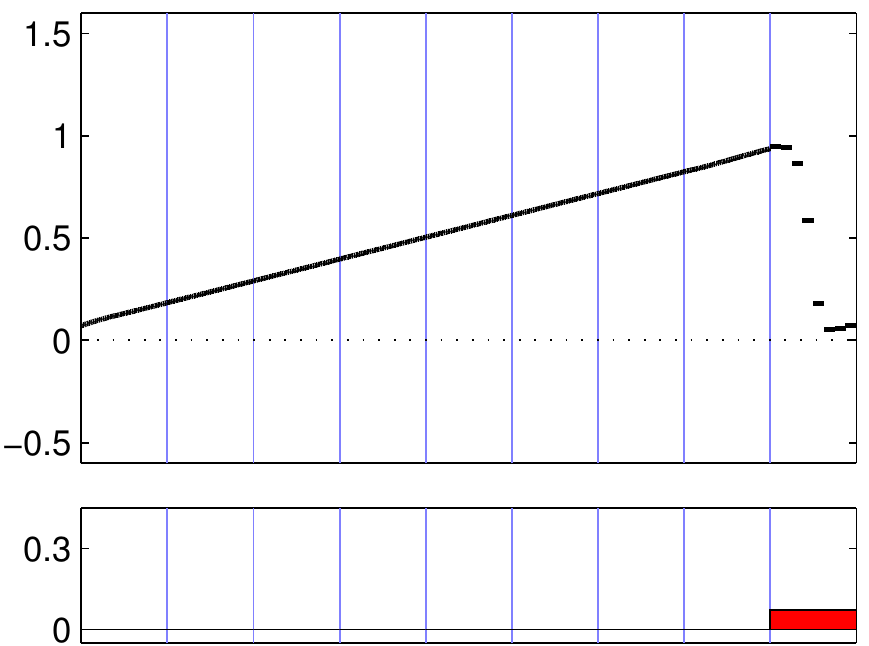} \\
    $t=0.55$ & $t=0.66$ & $t=0.77$ & $t=0.88$ \\
  \end{tabular}
  \caption{Evolution of the solution to the inviscid Burgers' with $p=4$, $n=8$ and $n_T=9$ elements at various times.}
  \label{fig:burgers}
\end{figure}

We now consider the inviscid Burgers' equation,  see Section \ref{ssec:GovEqns} {\tt [2]}, on $\Omega=(0,1)$ with smooth initial condition $u_0(x) = 1/2+\sin(2\pi x)$. We illustrate in Figure~\ref{fig:burgers} the solution and the indicator $\gamma$ with $p=4$, $n=8$ on $n_T=9$ elements at various times. The time-step is set to $\Delta t=10^{-3}$.
We observe that the indicator starts to be activated when the strong gradient appears even before the shock is formed since a polynomial representation with $p=4$ of the solution would yield oscillations. 
Second, it also can be seen that the indicator can deal with the moving shock and that the solution is recovered in polynomial representation after the shock left an element.
We also note that the standard discontinuous Galerkin method with high order elements will not lead to an approximation since over- and undershoots are amplified and the simulation aborted.

\subsection{Transonic Quasi-1D flow through nozzle}

We now shed our attention to a nozzle flow-problem under quasi-1D
assumptions. The resulting equations are similar to the one-dimensional Euler
equations, see Section \ref{ssec:GovEqns} {\tt [3]}, with minor modifications.
For the variable-area $A(x)$, this leads to the equations
\begin{equation} \label{eqn:nozzle-eqn}
 \begin{aligned}
  \pder{}{x}(A\rho u) &= 0, \\
  \pder{}{x}(A[\rho u^2 + p]) &= \frac{p}{A}\pder{A}{x}, \\
  \pder{}{x}(A[\rho E+p]u) &= 0
 \end{aligned}
\end{equation}
where $\rho$ is the fluid density, $u$ is the fluid velocity, and $p$ is the
thermodynamic pressure. The total energy is given by $ \rho E = \rho e + \rho
u^2 / 2$, and the pressure is related to $\rho E$ by the equation of state for
a perfect gas, $p = (\gamma_{\mathsf{a}}-1)\left(\rho E - \rho u^2 / 2 \right)$, where the
ratio of specific heats $\gamma_{\mathsf{a}} = 1.4$. The Mach number is the
ratio between the speed of the flow and the speed of sound $M=u/c$, where
$c=\sqrt{\gamma_ap/\rho}$.

We solve on the domain $\Omega=(0,1)$,
and prescribe the nozzle area by
\begin{equation}
 A(x) =
 \begin{cases}
  1 - (1-T_{\mathsf{n}})\cos\left(\pi\frac{x-0.5}{0.8}\right)^2 \quad &\text{for }0.1\le x \le 0.9 \\
  1 \quad &\text{otherwise},
 \end{cases}
\end{equation}
where $T_{\mathsf{n}} = 0.8$ is the height of the nozzle throat. At $x=0$ we impose the
farfield conditions $\rho=\rho_i = 1.0$, $u=u_i = 1.0$, $M=M_i = 0.40$ weakly, and at
$x=1$ we impose the conditions $\rho=\rho_o = 1.0$, $u=u_o = 1.0$, $M=M_o = 0.45$.

\begin{figure}
  \centering
  \begin{tabular}{cccc}
    \includegraphics[width=.23\textwidth]{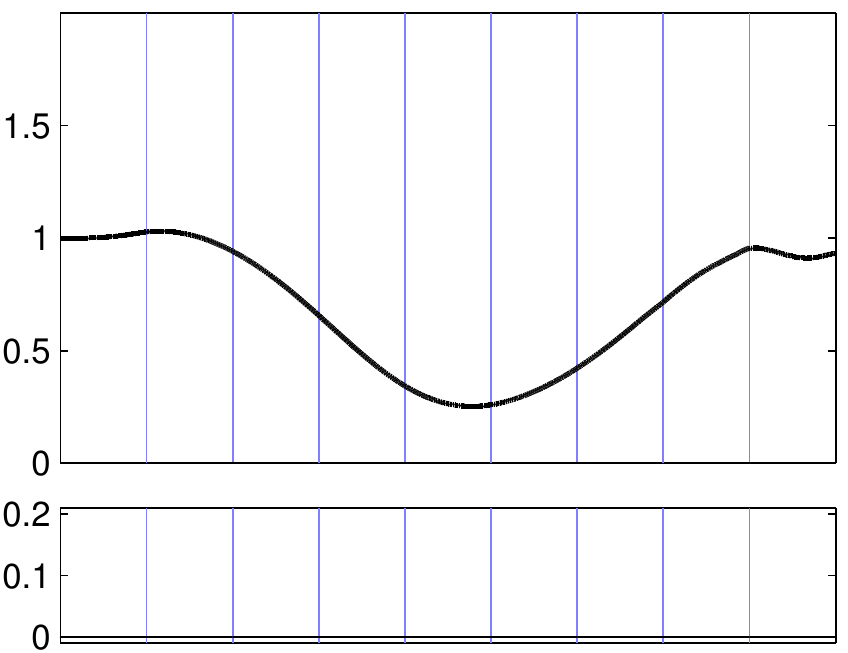} &
    \includegraphics[width=.23\textwidth]{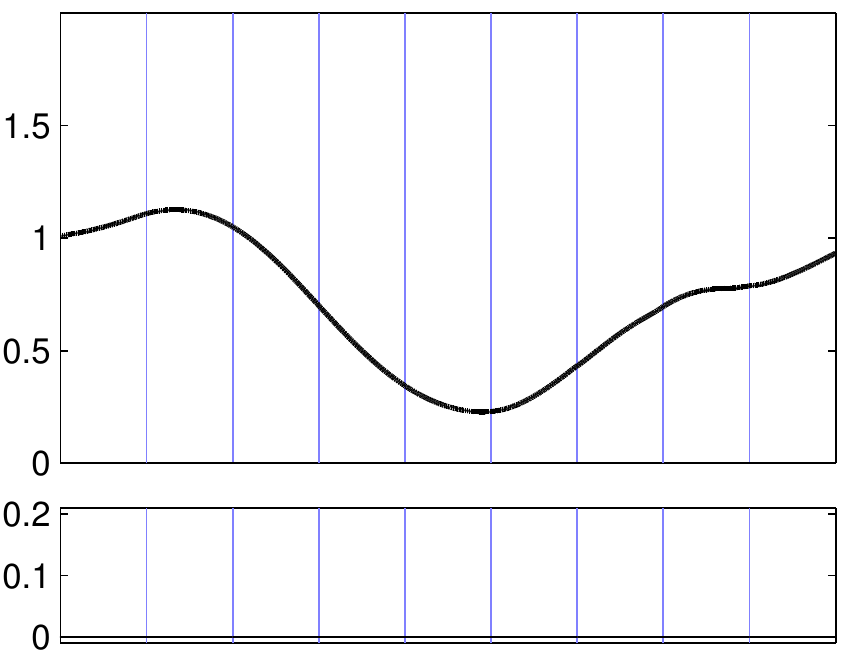} &
    \includegraphics[width=.23\textwidth]{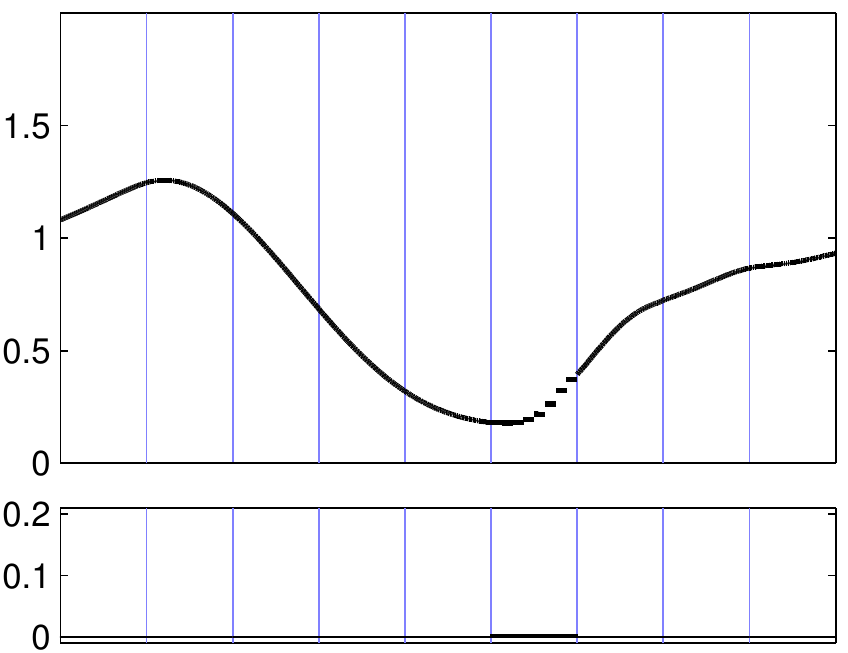} &
    \includegraphics[width=.23\textwidth]{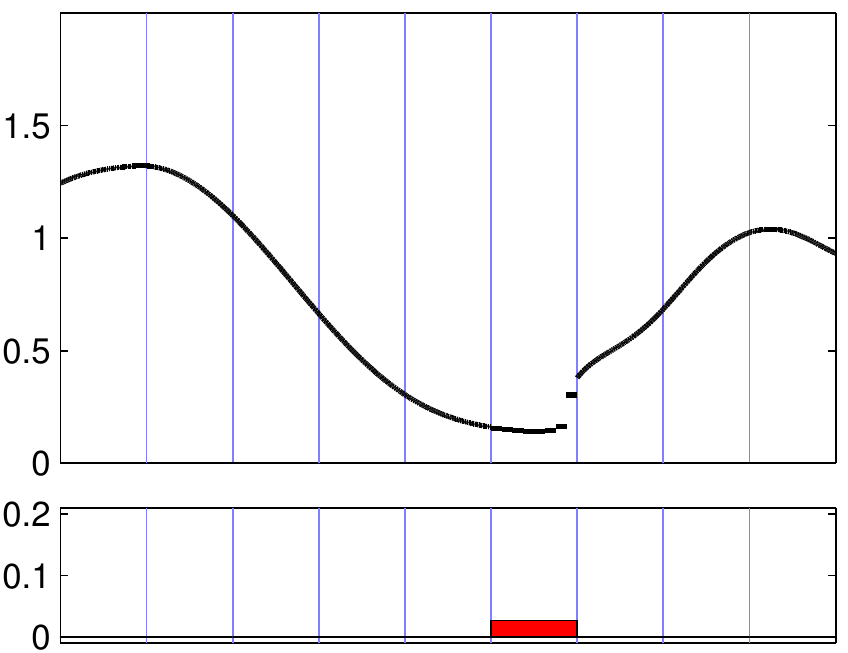} \\
    $t=0.05$ & $t=0.10$ & $t=0.15$ & $t=0.20$ \\ \\
    \includegraphics[width=.23\textwidth]{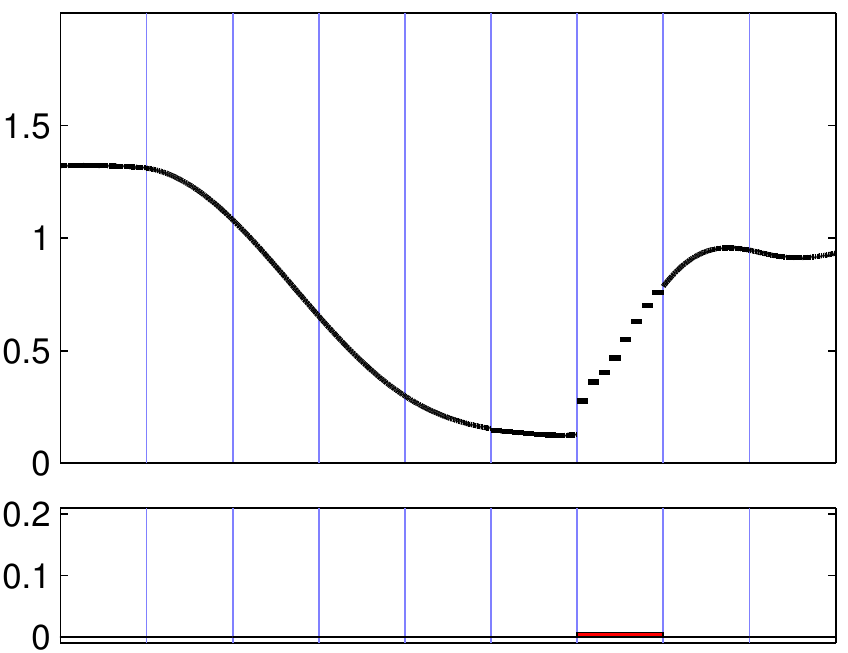} &
    \includegraphics[width=.23\textwidth]{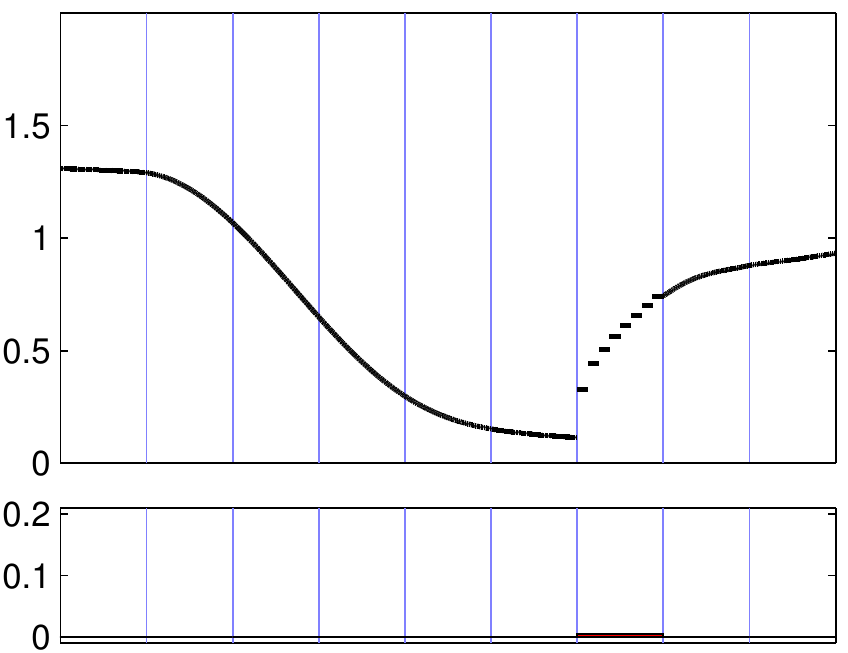} &
    \includegraphics[width=.23\textwidth]{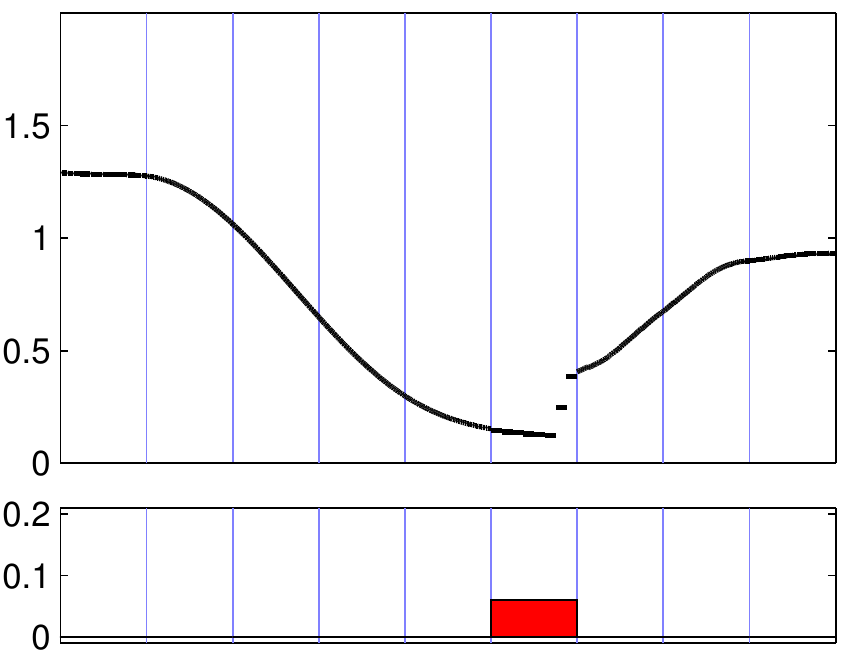} &
    \includegraphics[width=.23\textwidth]{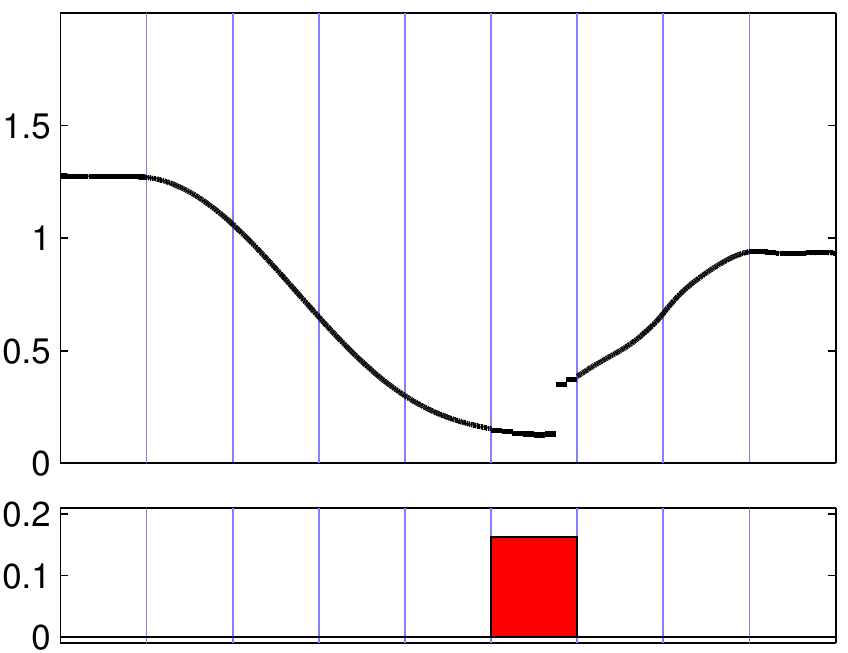} \\
    $t=0.25$ & $t=0.30$ & $t=0.35$ & $t=0.40$ \\
  \end{tabular}
  \caption{Transonic flow through a quasi-1D nozzle. The plots show
    the weighted density $A\rho(x)$ at various times. Discretized using
    $p=4$ and $n=8$. Note the sharp resolution of the shock on the
    sub-grid.}
  \label{fig:nozzle}
\end{figure}

This problem is a challenging problem for high-order methods since
small densities occur and thus small undershoots lead to very large relative
undershoots. 
Figure \ref{fig:nozzle} presents the first solution component
$A\rho$ at various times with $p=4$, $n=8$ on $n_T=9$ elements. The time-step
is set to $\Delta t=2\cdot 10^{-4}$. We observe that the steady-state is
reached and that the indicator is activated in the element containing the
shock.

\subsection{The Shu-Osher Shock Tube Problem}

Next we model the Shu-Osher problem \cite{shuosher}, which consists of a shock
front moving inside a one-dimensional inviscid flow with artificial density
fluctuations.

The governing equations are the 1D compressible Euler equations for an ideal
gas with a constant ratio of specific heats equal to $\gamma = 1.4$. The flow
domain is $\Omega = (-5,5)$ and the time domain is $(0, T)$ with $T =
1.78$. The initial condition for the flow is
\begin{align}
  (\rho,u,p) &=
  \begin{cases}
    (3.857143, 2.629369, 10.3333) & x<-4 \\
    (1 + 0.2 \sin(5x), 0, 1) & x\ge -4
  \end{cases}
\end{align}
with a supersonic inlet at $x = -5$ that prescribes the density, velocity, and
pressure $(\rho, u, p) = (3.857143, 2.629369, 10.3333)$ and a solid wall at $x
= 5$. This problem corresponds to a Mach $M = 3$ shock moving into a field with
a small density (or entropy) disturbance.

We solve the problem for polynomial degrees $p=1,2,3$ on two different
uniform grids, with $n_T=64$ and $n_T=256$ elements, respectively. In all
cases, we set the sub-grid size to $n=p+2$. The timestep is chosen small
enough that temporal errors are negligble.

\begin{figure}
  \centering
  \begin{tabular}{cc}
    \includegraphics[width=.46\textwidth]{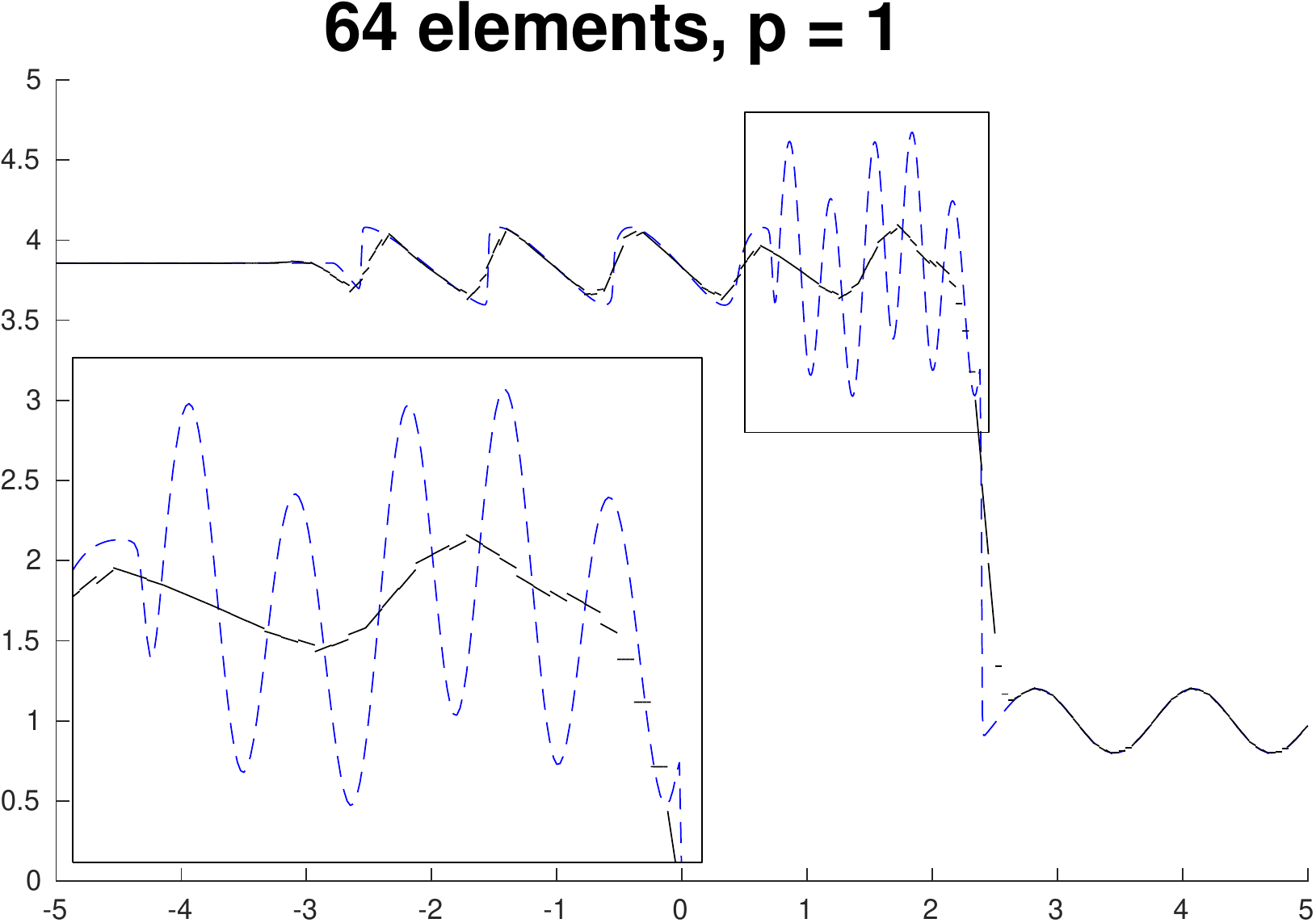} &
    \includegraphics[width=.46\textwidth]{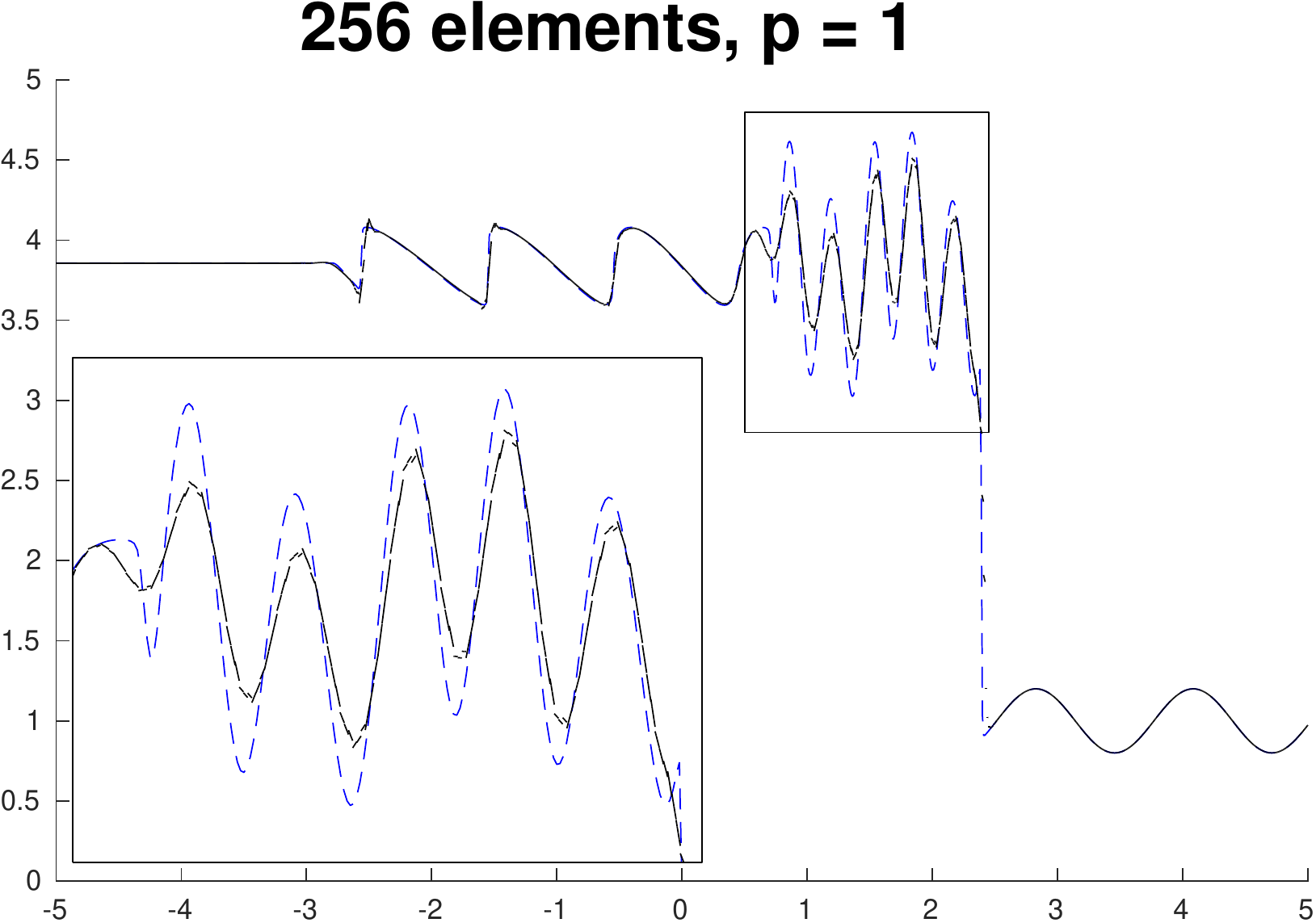} \\ \ \\
    \includegraphics[width=.46\textwidth]{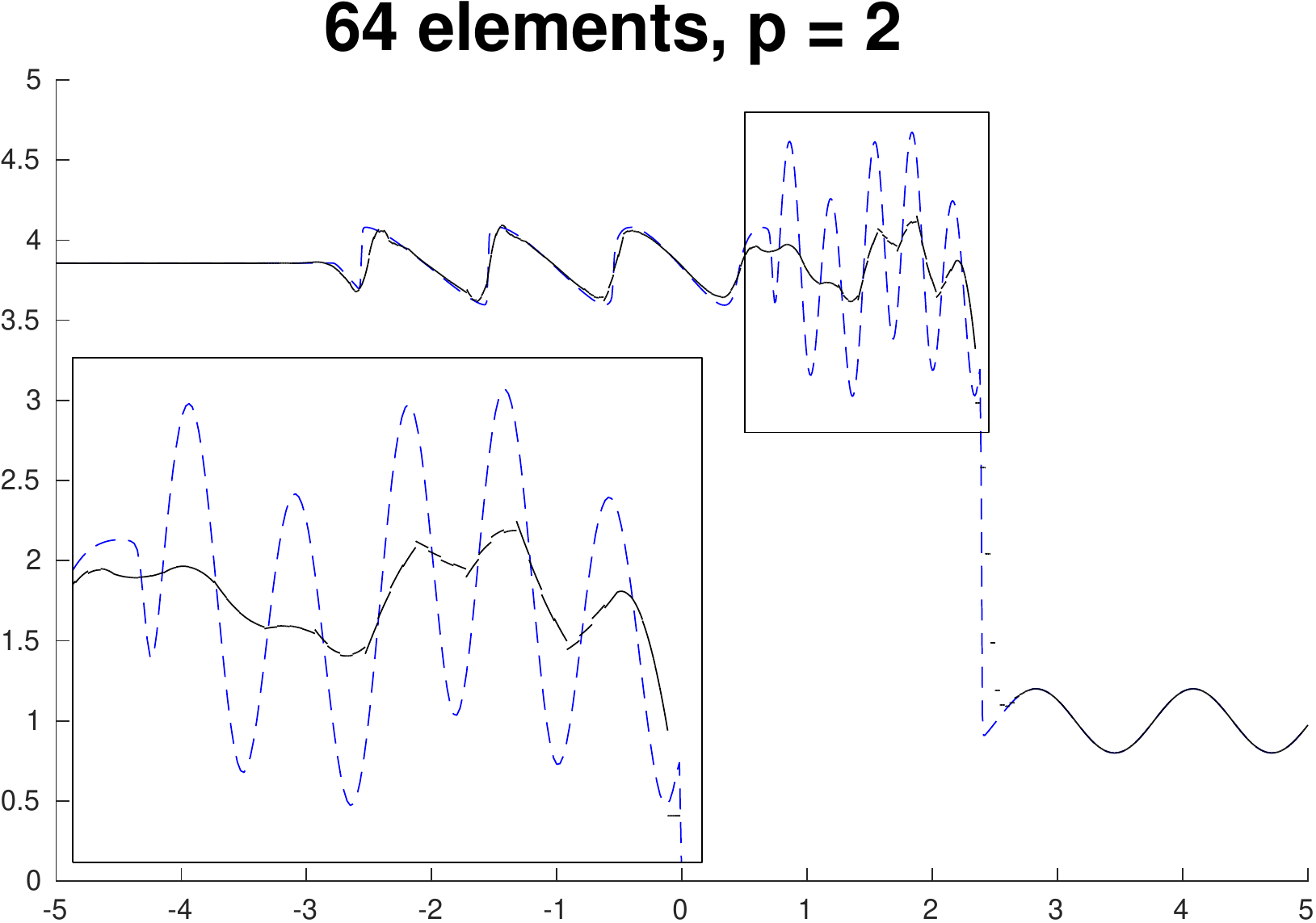} &
    \includegraphics[width=.46\textwidth]{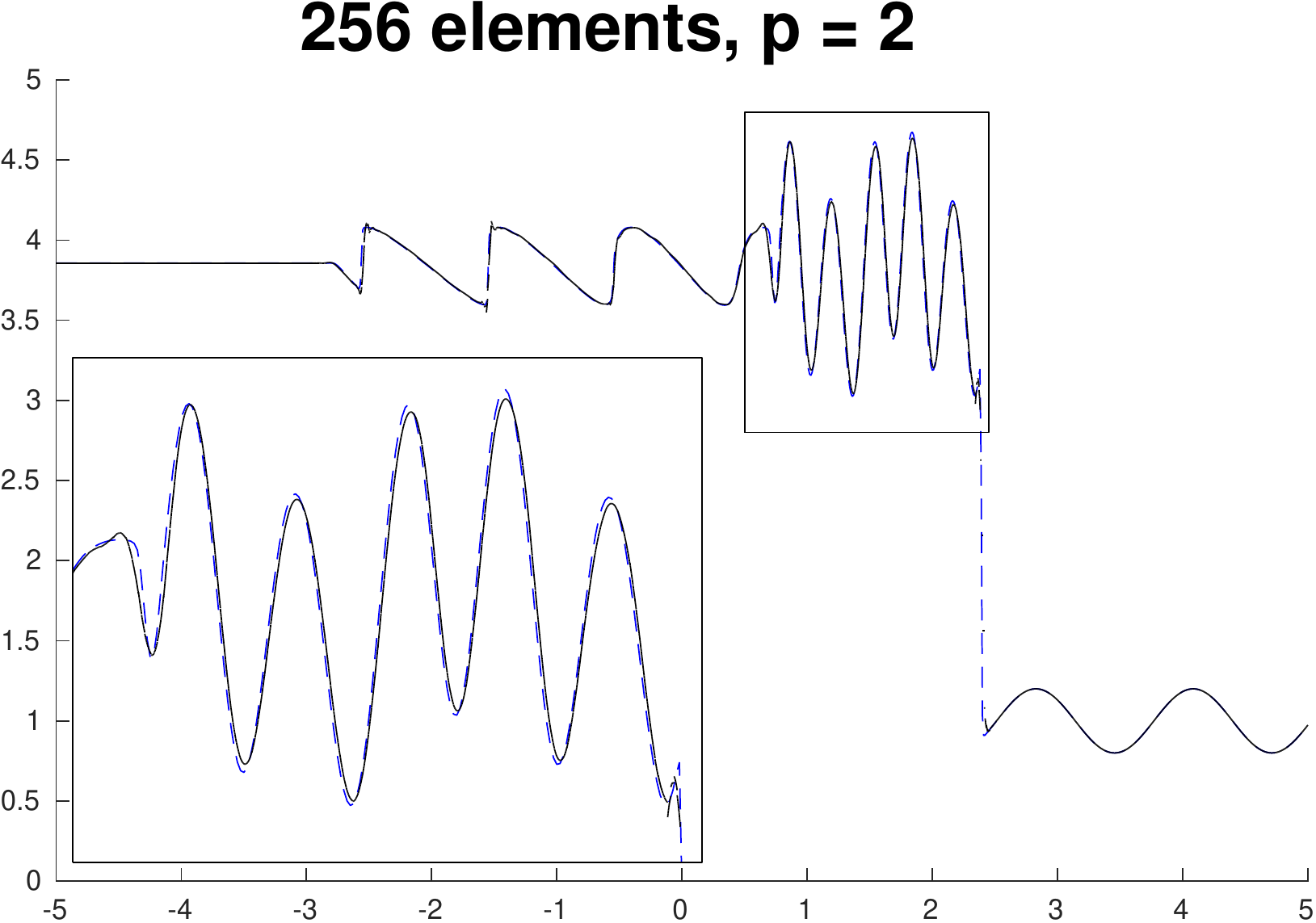} \\ \ \\
    \includegraphics[width=.46\textwidth]{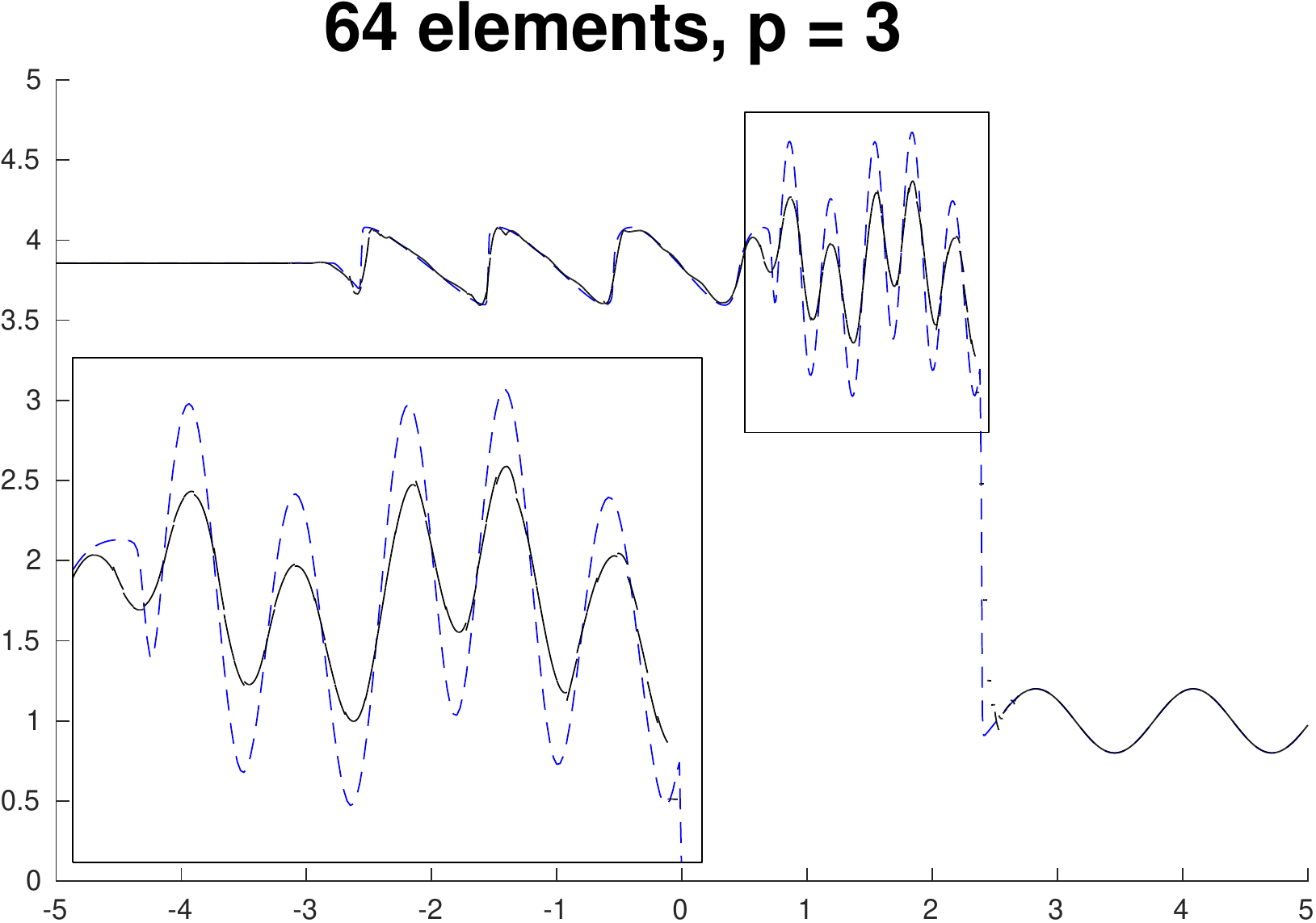} &
    \includegraphics[width=.46\textwidth]{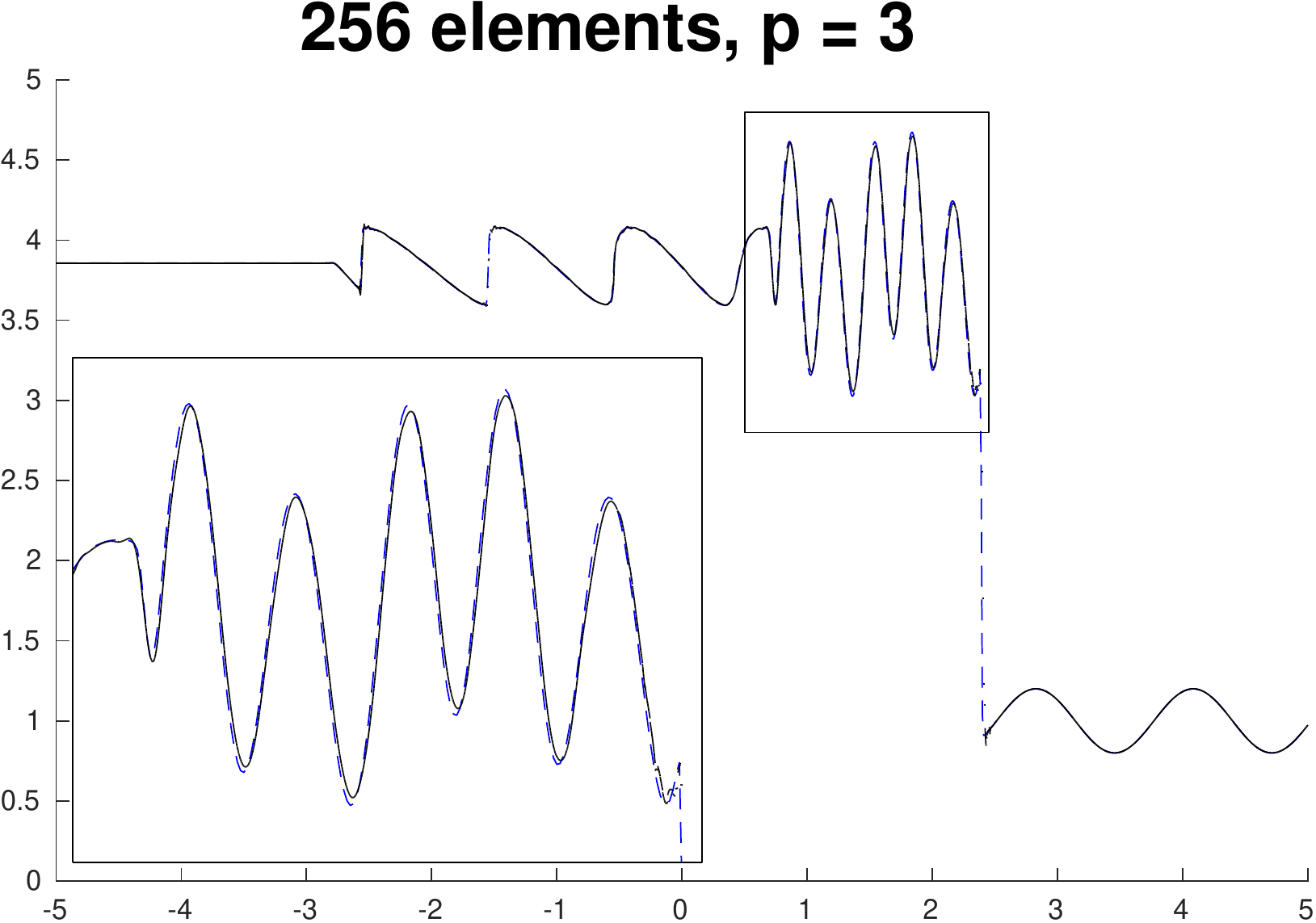} \\
  \end{tabular}
  \caption{The Shu-Osher Shock Tube problem. The plots show
    the density $\rho(x)$ at the final time $t=1.8$. Discretized using
    $p=1,2,3$ and $n=p+2$. The reference solutions is shown in blue
    dashed line.}
  \label{fig:shuosher}
\end{figure}

The resulting density fields are shown in Figure~\ref{fig:shuosher},
together with a reference solution (computed using a WENO scheme
on a very fine grid). We observe that the method generates stable
solutions for all cases, and that the accuracy of the scheme
increases substantially for higher polynomial degrees $p$.

\subsection{Euler Vortex}

Our first 2D example is the inviscid flow of a compressible vortex in a
rectangular domain \cite{erlebacher97vortex}. The vortex is initially centered
at $(x_0,y_0)$ and is moving with the free-stream at an angle $\theta$ with
respect to the $x$-axis. The analytic solution at $(x,y,t)$ is given by
\begin{align}
u &= u_\infty \left(\cos\theta - \frac{\epsilon ((y-y_0)-\bar{v}t)}{2\pi r_c}
                    \exp(f/2) \right), &
\rho &= \rho_\infty \left(1 - \frac{\epsilon^2(\gamma-1)M_\infty^2}{8\pi^2}
                    \exp (f) \right)^\frac{1}{\gamma-1}, \\
v &= u_\infty \left(\sin\theta + \frac{\epsilon ((x-x_0)-\bar{u}t)}{2\pi r_c}
                    \exp (f/2) \right), & 
p    &= p_\infty \left(1 - \frac{\epsilon^2(\gamma-1)M_\infty^2}{8\pi^2}
                    \exp (f) \right)^\frac{\gamma}{\gamma-1},
\end{align}
where $f(x,y,t) =
(1-((x-x_0)-\bar{u}t)^2-((y-y_0)-\bar{v}t)^2)/r_c^2$, $M_\infty$ is
the Mach number, $\gamma=c_p/c_v=1.4$, and $u_\infty$, $p_\infty$,
$\rho_\infty$ are free-stream velocity, pressure, and density. The
Cartesian components of the free-stream velocity are
$\bar{u}=u_\infty\cos\theta$ and $\bar{v}=u_\infty\sin\theta$. The
parameter $\epsilon$ measures the strength of the vortex and $r_c$ is
its size.

We use a domain of size 20-by-15, with the vortex initially centered at
$(x_0,y_0)=(5,5)$ with respect to the lower-left corner. The Mach number is
$M_\infty=0.5$, the angle $\theta=\arctan 1/2$, and the vortex has the
parameters $\epsilon=10$ and $r_c=1.5$. We use characteristic boundary
conditions and integrate until time $t_0=\sqrt{10^2+5^2}/10$, when the vortex
has moved a relative distance of $(1,1/2)$.

\begin{figure}[t]
  \begin{center}
    \includegraphics[width=.7\textwidth]{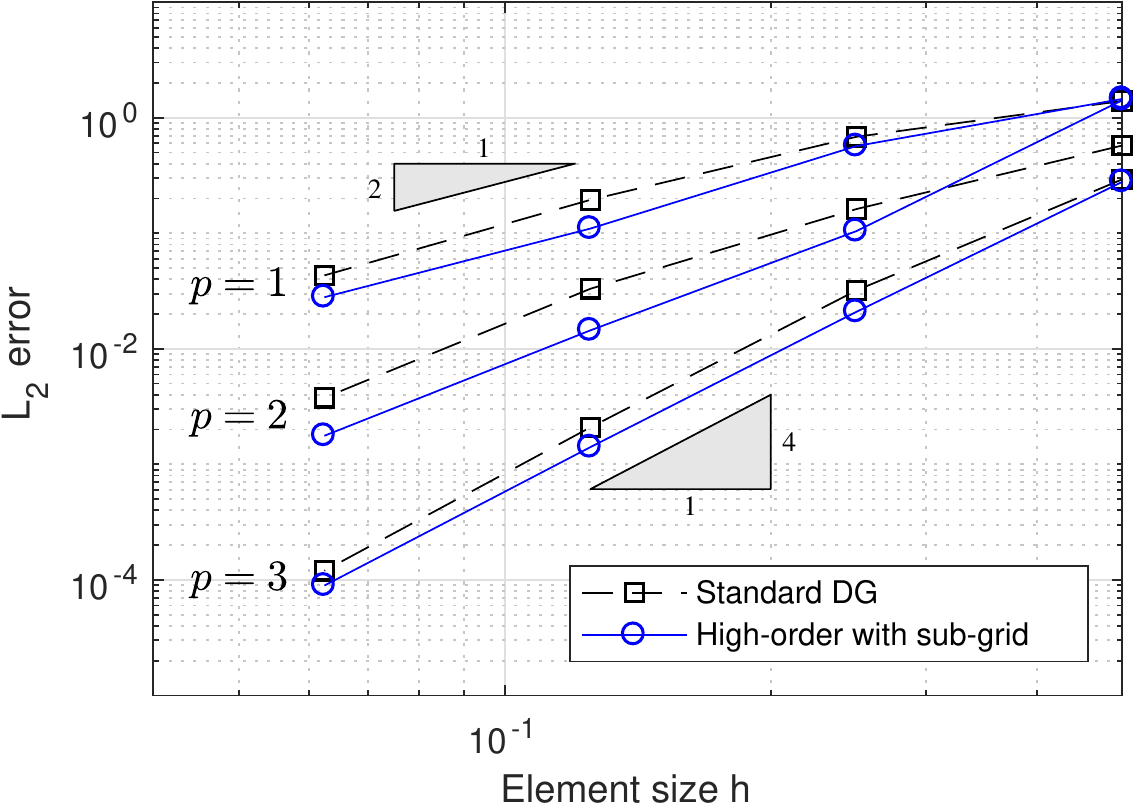}
  \end{center}
  \caption{Convergence test for an Euler vortex test problem using a standard
    DG method and our proposed mixed high- and low-order scheme.  The results
    show optimal order of convergence $\mathcal{O}(h^{p+1})$ for both methods.
    In general, the mixed scheme is more accurate than standard DG, since its
    approximation space is richer.}
  \label{eulervortex_convergence}
\end{figure}

In Figure~\ref{eulervortex_convergence}, we graph the $L_2$-errors for all
simulation cases, both for the standard DG method and our proposed method with
a range of mesh sizes $h$ and polynomial degrees $p$. The results clearly show
the optimal order of convergence $\mathcal{O}(h^{p+1})$ for both methods.
The new method is slightly more accurate, which is expected
since its approximation space is richer but contains the standard DG polynomials.

\subsection{The Woodward-Colella forward facing step}

We finally apply our method to the forward facing step problem of Woodward and
Colella \cite{colella84shocks}, in the context of Euler's equations in two
dimensions. The freestream Mach number is 3, and we discretize the high-order
space using polynomials of degree $p=3$, and create a sub-grid by refining
each triangle into $n=4^2=16$ sub-cells uniformly. We integrate in time using a
stepsize of $\Delta t=2\cdot 10^{-4}$ until a final time $5.0$. We only use a
first-order accurate IMEX scheme in time, since our goal is to compute a steady-state
solution which can be used to assess the shock capturing capabilities.
The solution never reaches a steady-state due to transient effects and
instabilities, but it serves as a good test case.

\begin{figure}
  \centering
   \includegraphics[width=.49\textwidth]{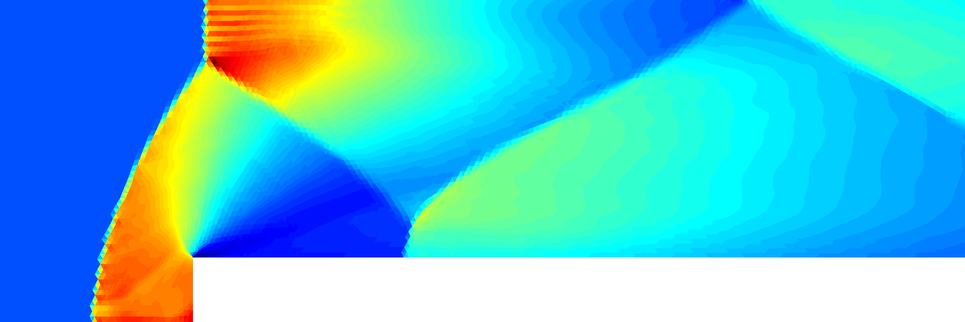} \hfill
   \includegraphics[width=.49\textwidth]{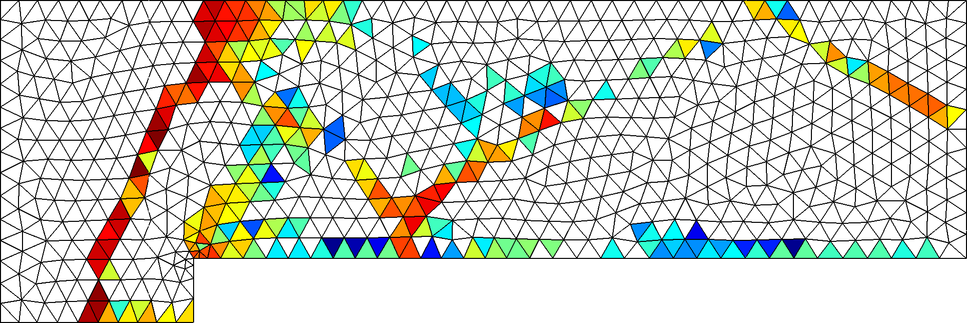} \\ \ \\
   \includegraphics[width=.49\textwidth]{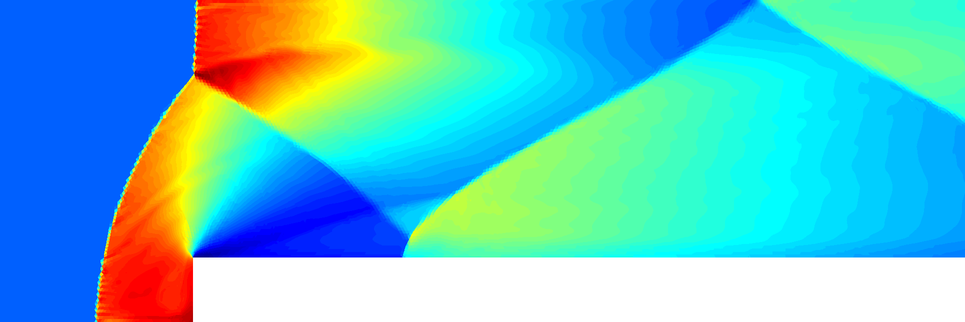} \hfill
   \includegraphics[width=.49\textwidth]{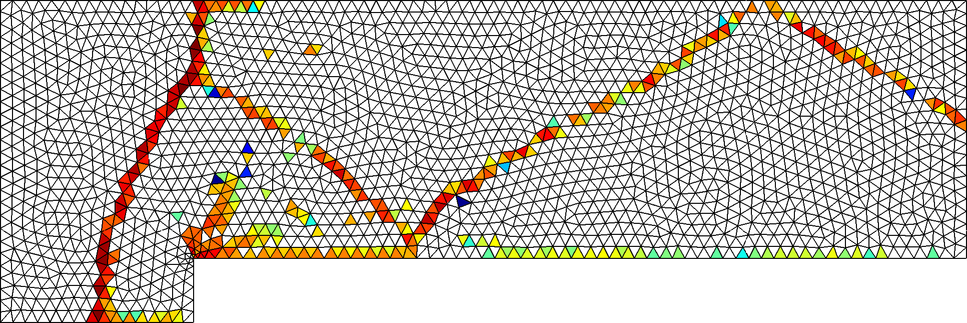} \\
  \caption{The Woodward-Colella forward facing step problem
    \cite{colella84shocks}. The problem is discretized on a coarse mesh (top)
    and on a finer mesh (bottom), using polynomial degrees $p=3$ and a
    sub-grid of $n=4^2=16$ sub-cells in each triangle. The plots show the density
    (left) and the sensor on the mesh (right).  In the sensor plots, gray
    color indicated no stabilization is applied. Note that most of the shocks
    are captured using only a single layer of low-order elements. The
    numerical oscillations are due to the low-order scheme (the so-called
    ``carbuncle effect'', \cite{Quirk1997,carbuncle2001}).}
  \label{fig:colella}
\end{figure}

Figure~\ref{fig:colella} shows the results on a coarse and on a finer mesh,
visualized by a density plot as well as the sensor and the mesh. First, we
note that our sensor is highly selective and that only a single element layer
of elements needs to be flagged for stabilization. The singularity at the
convex corner causes some problems, as the oscillations due to
under-resolution are convected downstream and results in an entire layer of
elements along the wall being flagged by the sensor. This could likely be
avoided by a finer resolution at the corner (as in \cite{persson13shock}),
different parameters for the sensor, or by special treatment of the corner as
in the original work \cite{colella84shocks}. We also note numerical
oscillations behind the shocks, in particular for the coarse mesh. This is due
to the low-order finite volume scheme and known as the ``carbuncle effect''
\cite{Quirk1997}, which can be remedied using a wide range of techniques
\cite{carbuncle2001}.

\subsection{Comparison with other methods}

Although it is difficult to draw conclusions about how our proposed high-order
with sub-grid scheme compares to other methods, we can still illustrate how the
solutions behave qualitatively and give some comments about the computational
cost of the method.

\begin{figure}
  \centering
  \begin{tabular}{cc}
    \includegraphics[width=.46\textwidth]{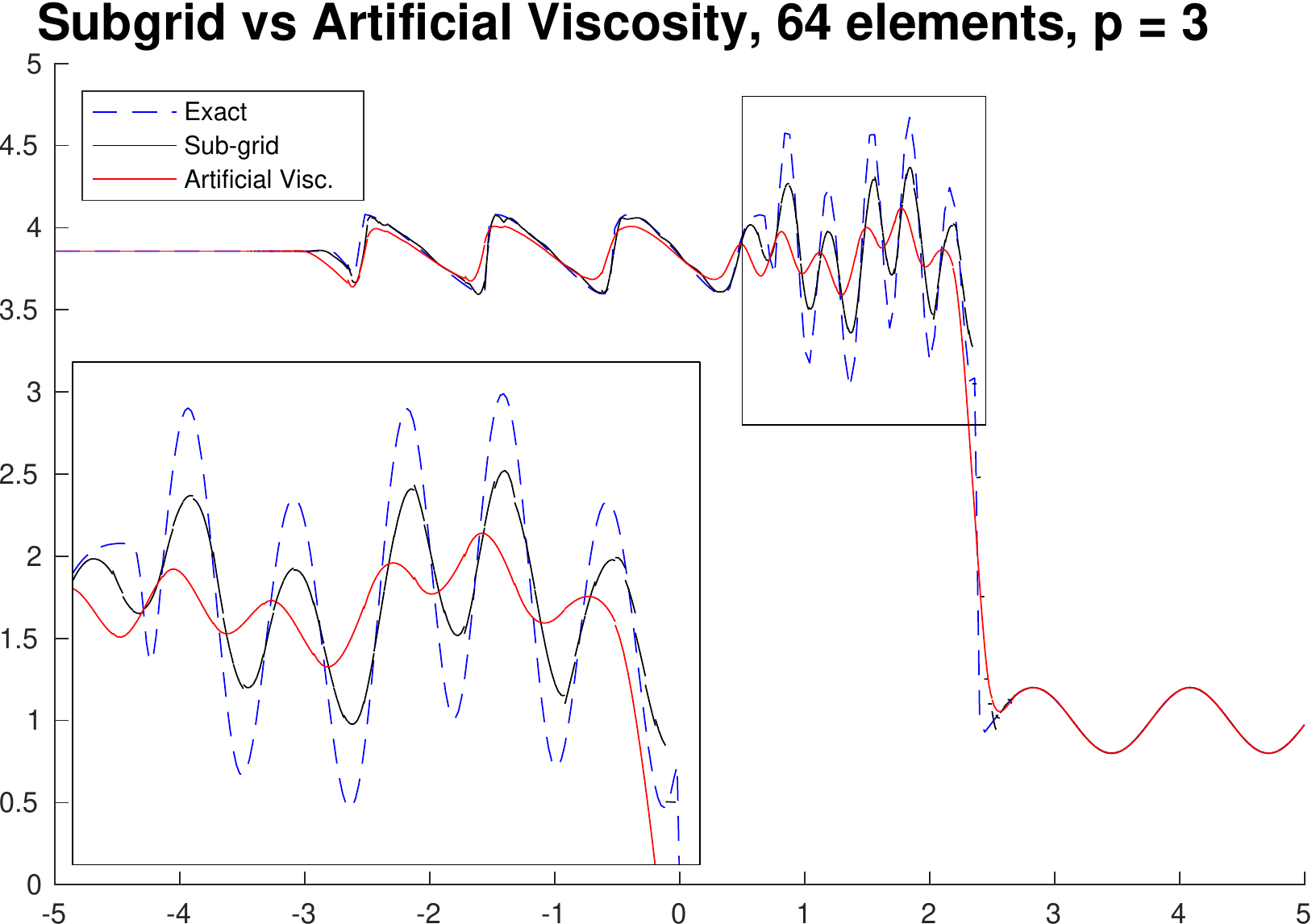} &
    \includegraphics[width=.46\textwidth]{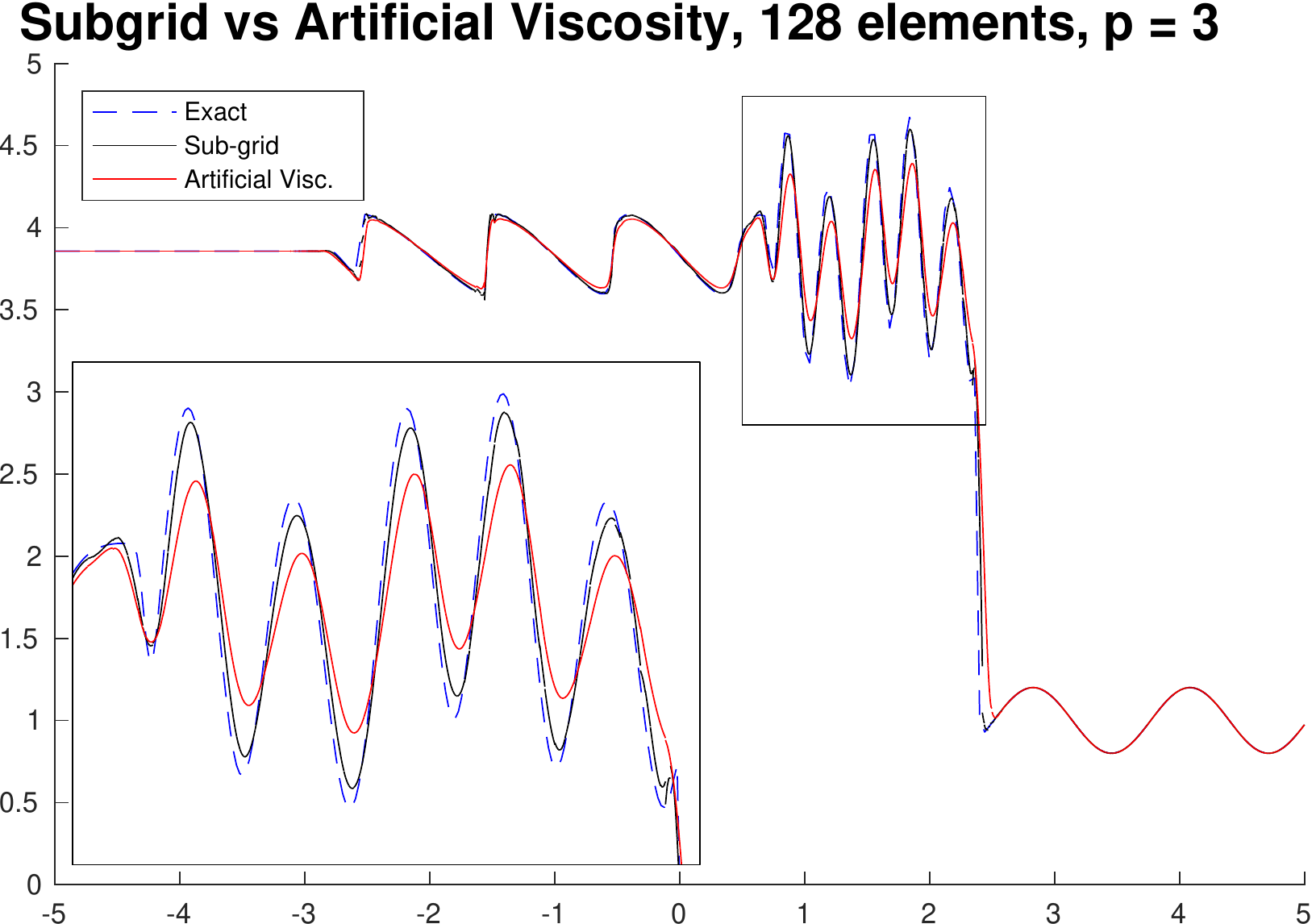}
  \end{tabular}
  \caption{Comparison between the high-order with sub-grid method and standard
    DG with artificial viscosity shock
    capturing~\cite{persson13shock,persson06shock}, for the Shu-Osher shock
    tube problem.}
  \label{fig:so1cmpav}
\end{figure}

In Figure~\ref{fig:so1cmpav}, we show the Shu-Osher shock tube problem
introduced earlier, using both the new sub-grid scheme and a standard DG scheme
with the artificial viscosity stabilization proposed
in~\cite{persson13shock,persson06shock}. The amount of viscosity is chosen as
small as possible to give stable solutions. Still, the sub-grid method gives
much less numerical dissipation as well as a sharper shock profile.

\begin{figure}
  \centering
  \begin{tabular}{cc}
    \includegraphics[width=.46\textwidth]{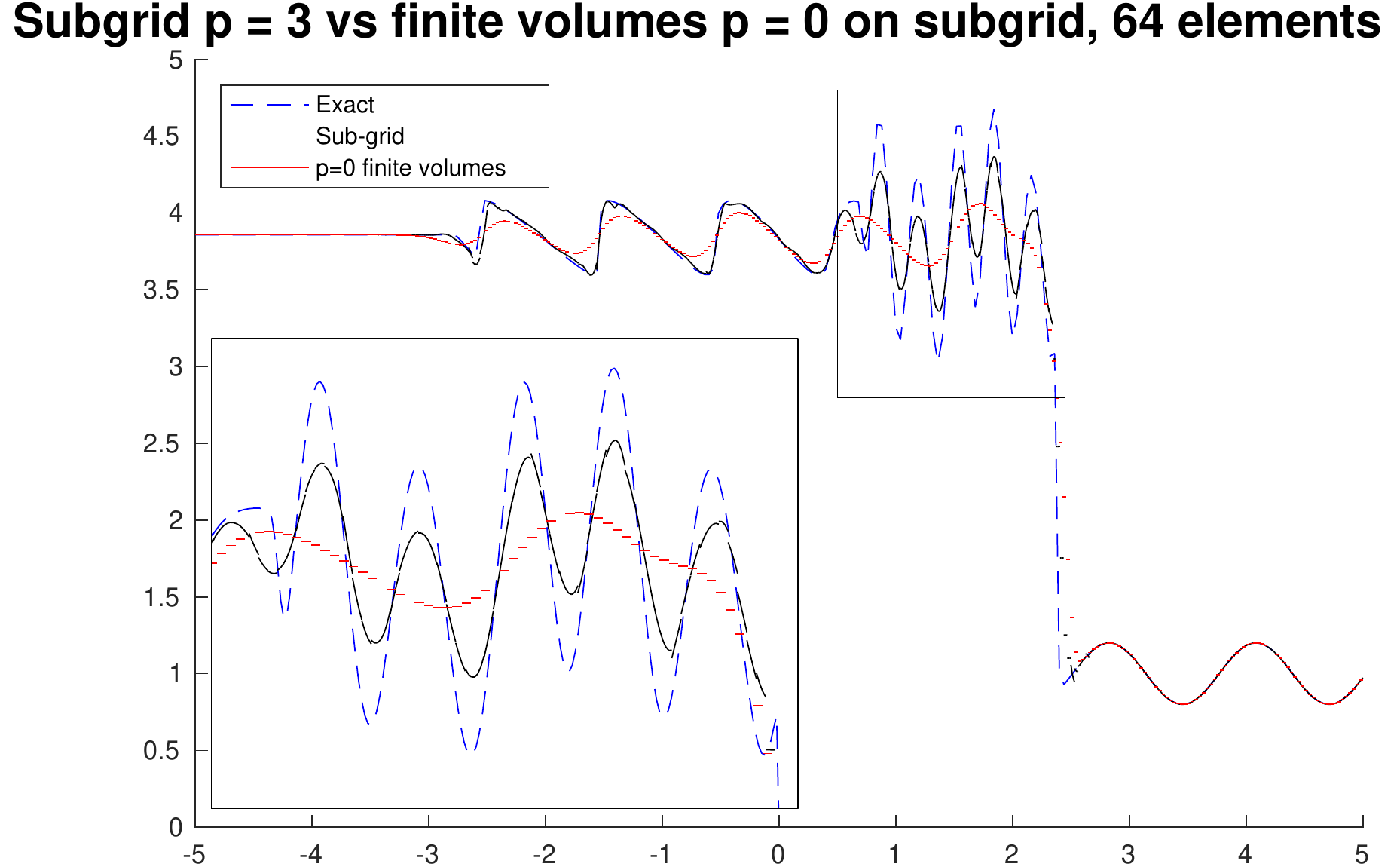} &
    \includegraphics[width=.46\textwidth]{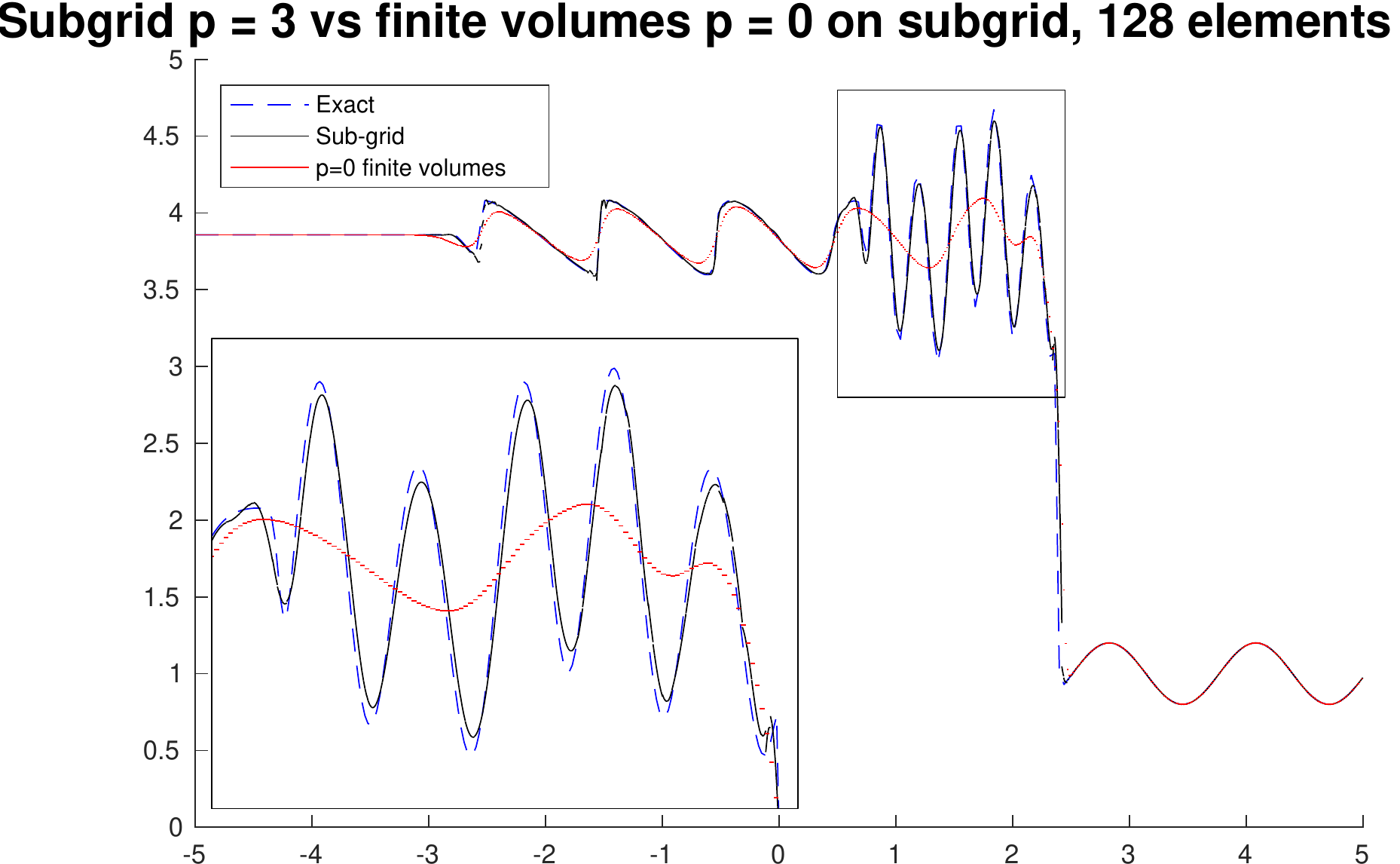}
  \end{tabular}
  \caption{Comparison between the high-order with sub-grid method and
    cell-centered 1st order finite volumes (standard DG with $p=0$) on the
    sub-grid, for the Shu-Osher shock tube problem.}
  \label{fig:so1cmpfvm}
\end{figure}

To illustrate the benefits of the high-order accuracy, we also compare the
with a standard cell-centered finite volume scheme in
Figure~\ref{fig:so1cmpfvm}. The grid is the same as the sub-grid used for the
new method, that is, it corresponds to setting $p=0$ in our sub-grid
method. The finite volume solution is highly inaccurate, in particular in the
oscillatory region behind the shock, showing the benefit of the high-order
polynomial modes.

Finally, we make some comments on the performance of our proposed scheme.  Our
implementation is only a prototype code, so wall-clock time comparisons with
highly optimized codes would not give a accurate view of the actual
computational cost. Instead, let us just comment on some of the main
considerations:

\begin{description}
\item[Degrees of freedom:] Clearly, the high-order with sub-grid method that we
  propose uses more degrees of freedom than a corresponding standard DG method
  on the same grid.  More precisely, it needs $n$ more degrees of freedom per
  solution component and element. Since we choose the number of sub-grid
  elements $n$ approximately the same as the number of high-order DG
  coefficients, this ends up being about twice the number of degrees of
  freedom.  However, note from the convergence plot in
  Figure~\ref{eulervortex_convergence} that these extra degrees of freedom
  also result in a more accurate solution.

\item[Computational cost of residual assembly:] Our implementation uses a
  standard DG code to compute all the integrals on the sub-grid, which is
  clearly inefficient since it assumes full polynomial degrees $p$ everywhere
  (i.e. even on the sub-grid).  This can easily be improved in a specialized
  code for the proposed sub-grid scheme.  However, an important decision is how
  to compute the integrals. The most straight-forward way is to use a
  high-order quadrature rule inside each sub-grid element and on the edges
  between them, which results in about $n$ times as many quadrature points as
  standard DG. Future work includes the development of specialized quadrature
  rules for our particular solution spaces.

\item[Adaptive assembly of elements with shocks:] Typically, the high-order
  with sub-grid shock capturing scheme is only active in a small fraction of
  the mesh elements. Therefore, an obvious way to improve the performance of
  the method is to adaptively only apply the new scheme in those elements.
  This might need a more efficient shock indicator than the one we propose
  here, for example the resolution based indicator in~\cite{persson06shock},
  to avoid having to assemble the sub-grid everywhere. This should in principle
  make the cost of our method comparable with standard DG for most practical
  problems.
\end{description}

\section{Conclusions}
In this article, we have introduced a new numerical scheme to solve non-linear conservation laws with a discontinuous Galerkin method based upon a non-standard (discontinuous) approximation space. This space includes on each element both high-order polynomials and piecewise constant functions on a sub-grid. 
The set-up of the method is very general so that it can be applied to unstructured meshes consisting of simplex elements in any dimension.

Since a Galerkin method is a projection method which minimizes the residual over the approximation space, we argue that the method intrinsically has a tendency to choose a good balance between using low- or high-order features of the approximation space. 
This, however, only reduces, but does not cure, the problem of overshooting which is fundamental for non-linear problems. 
The particular structure of the approximation space allows to define a sensor with some nice properties which allows to define a local penalty parameter that suppresses the high order modes whenever the sensor is activated in a smooth way. The method also allows for recovering a polynomial representation once a shock has quit an element, however with an accuracy that is reduced due to the temporary low-order representation of the solution. 

Several numerical tests illustrate the characteristics of the method. From our prototypical implementation, it is difficult to guess how expensive the method really is. We judge it however satisfying that we obtain a qualitatively correct solution in regard that a typical high-order method crashes without any sort of shock-capturing and that the high-order modes are only suppressed locally where needed.

\section{Acknowledgement}
Financial support from the France-Berkeley Fund under the project-name ``Sub-cell Limiting for Shock-capturing in High-Order Discontinuous Galerkin Methods'' is acknowledged.

\bibliographystyle{plain}
\bibliography{Biblio}

\begin{thebibliography}{10}

\bibitem{barter08}
Garrett Barter.
\newblock {\em Shock capturing with {PDE}-based artificial viscosity for an
  adaptive, higher-order discontinuous {G}alerkin finite element method}.
\newblock PhD thesis, M.I.T., June 2008.

\bibitem{bassirebay95shocks}
Francesco Bassi and S.~Rebay.
\newblock Accurate {2D} {E}uler computations by means of a high order
  discontinuous finite element method.
\newblock In {\em Fourteenth International Conference on Numerical Methods in
  Fluid Dynamics}, pages 234--240. Springer, 1995.

\bibitem{BauOden}
C.~E. Baumann and J.~T. Oden.
\newblock A discontinuous {$hp$} finite element method for the {E}uler and
  {N}avier-{S}tokes equations.
\newblock {\em Int. J. Numer. Methods Fluids}, 31(1):79--95, 1999.
\newblock Tenth International Conference on Finite Elements in Fluids (Tucson,
  AZ, 1998).

\bibitem{lele2009viscosity}
Ankit Bhagatwala and Sanjiva Lele.
\newblock A modified artificial viscosity approach for compressible turbulence
  simulations.
\newblock {\em Journal of Computational Physics}, 228(14):4965--4969, 2009.

\bibitem{gassner2019shocks}
Marvin Bohm, Sven Schermeng, Andrew~R. Winters, Gregor~J. Gassner, and
  Gustaaf~B. Jacobs.
\newblock Multi-element {SIAC} filter for shock capturing applied to high-order
  discontinuous {G}alerkin spectral element methods.
\newblock {\em J. Sci. Comput.}, 81(2):820--844, 2019.

\bibitem{burbeau01limiter}
A.~Burbeau, P.~Sagaut, and Ch.-H. Bruneau.
\newblock A problem-independent limiter for high-order {R}unge-{K}utta
  discontinuous {G}alerkin methods.
\newblock {\em Journal of Computational Physics}, 169(1):111--150, 2001.

\bibitem{BookQuarteroni}
C.~Canuto, M.~Y. Hussaini, A.~Quarteroni, and T.~A. Zang.
\newblock {\em Spectral methods}.
\newblock Scientific Computation. Springer-Verlag, Berlin, 2006.
\newblock Fundamentals in single domains.

\bibitem{cockburn01rkdg}
B.~Cockburn and C.-W. Shu.
\newblock Runge-{K}utta discontinuous {G}alerkin methods for
  convection-dominated problems.
\newblock {\em J. Sci. Comput.}, 16(3):173--261, 2001.

\bibitem{dervieux03adaptation}
Alain Dervieux, David Leservoisier, Paul-Louis George, and Yves Coudi{\`e}re.
\newblock About theoretical and practical impact of mesh adaptation on
  approximation of functions and {PDE} solutions.
\newblock {\em International Journal for Numerical Methods in Fluids},
  43(5):507--516, 2003.
\newblock ECCOMAS Computational Fluid Dynamics Conference, Part I (Swansea,
  2001).

\bibitem{dumbser2016simple}
Michael Dumbser and Rapha{\"e}l Loub{\`e}re.
\newblock A simple robust and accurate a posteriori sub-cell finite volume
  limiter for the discontinuous galerkin method on unstructured meshes.
\newblock {\em Journal of Computational Physics}, 319:163--199, 2016.

\bibitem{dumbser2014subcell}
Michael Dumbser, Olindo Zanotti, Rapha\"{e}l Loub\`ere, and Steven Diot.
\newblock {\it {A} posteriori} subcell limiting of the discontinuous {G}alerkin
  finite element method for hyperbolic conservation laws.
\newblock {\em J. Comput. Phys.}, 278:47--75, 2014.

\bibitem{erlebacher97vortex}
G.~Erlebacher, M.~Y. Hussaini, and C.-W. Shu.
\newblock Interaction of a shock with a longitudinal vortex.
\newblock {\em J. Fluid Mech.}, 337:129--153, 1997.

\bibitem{eno1}
Ami Harten, Bj\"orn Engquist, Stanley Osher, and Sukumar~R. Chakravarthy.
\newblock Uniformly high-order accurate essentially nonoscillatory schemes.
  {III}.
\newblock {\em Journal of Computational Physics}, 71(2):231--303, 1987.

\bibitem{hartmann02shock}
Ralf Hartmann and Paul Houston.
\newblock Adaptive discontinuous {G}alerkin finite element methods for the
  compressible {E}uler equations.
\newblock {\em Journal of Computational Physics}, 183(2):508--532, 2002.

\bibitem{huerta2012shockDG}
A.~Huerta, E.~Casoni, and J.~Peraire.
\newblock {A simple shock-capturing technique for high-order discontinuous
  Galerkin methods}.
\newblock {\em International Journal for Numerical Methods in Fluids},
  69(10):1614--1632, 2012.

\bibitem{jameson81jst}
Antony Jameson, W.~Schmidt, and Eli Turkel.
\newblock Numerical solution of the {E}uler equations by finite volume methods
  using {R}unge {K}utta time stepping schemes.
\newblock In {\em 14th Fluid and Plasma Dynamics Conference}, 1981.

\bibitem{weno2}
Guang-Shan Jiang and Chi-Wang Shu.
\newblock Efficient implementation of weighted {ENO} schemes.
\newblock {\em Journal of Computational Physics}, 126(1):202--228, 1996.

\bibitem{klockner2011shock}
Andreas Kl\"ockner, Tim Warburton, and Jan~S. Hesthaven.
\newblock Viscous shock capturing in a time-explicit discontinuous {G}alerkin
  method.
\newblock {\em Mathematical Modelling of Natural Phenomena}, 6(3):57--83, 2011.
\newblock arXiv:1102.3190.

\bibitem{weno1}
Xu-Dong Liu, Stanley Osher, and Tony Chan.
\newblock Weighted essentially non-oscillatory schemes.
\newblock {\em Journal of Computational Physics}, 115(1):200--212, 1994.

\bibitem{moro14thesis}
David Moro-Ludena.
\newblock {\em {An adaptive high order Reynolds-averaged Navier-Stokes solver
  with transition prediction}}.
\newblock PhD thesis, M.I.T., December 2014.

\bibitem{carbuncle2001}
Maurizio Pandolfi and Domenic D'Ambrosio.
\newblock Numerical instabilities in upwind methods: analysis and cures for the
  ``carbuncle'' phenomenon.
\newblock {\em J. Comput. Phys.}, 166(2):271--301, 2001.

\bibitem{persson06shock}
P.-O. Persson and J.~Peraire.
\newblock Sub-cell shock capturing for discontinuous {G}alerkin methods.
\newblock In {\em 44th AIAA Aerospace Sciences Meeting and Exhibit, Reno,
  Nevada}, 2006.
\newblock AIAA-2006-0112.

\bibitem{persson13shock}
Per-Olof Persson.
\newblock Shock capturing for high-order discontinuous {G}alerkin simulation of
  transient flow problems.
\newblock In {\em 21st AIAA Computational Fluid Dynamics Conference, San Diego,
  CA}, Jun 2013.
\newblock AIAA-2013-3061.

\bibitem{Quirk1997}
James~J. Quirk.
\newblock A contribution to the great {R}iemann solver debate.
\newblock In {\em Upwind and High-Resolution Schemes}, pages 550--569. Springer
  Berlin Heidelberg, Berlin, Heidelberg, 1997.

\bibitem{BookSchwab}
Ch. Schwab.
\newblock {\em {$p$}- and {$hp$}-finite element methods}.
\newblock Numerical Mathematics and Scientific Computation. The Clarendon
  Press, Oxford University Press, New York, 1998.
\newblock Theory and applications in solid and fluid mechanics.

\bibitem{shuosher}
Chi-Wang Shu and Stanley Osher.
\newblock Efficient implementation of essentially nonoscillatory
  shock-capturing schemes. {II}.
\newblock {\em J. Comput. Phys.}, 83(1):32--78, 1989.

\bibitem{sonntag2017efficient}
Matthias Sonntag and Claus-Dieter Munz.
\newblock Efficient parallelization of a shock capturing for discontinuous
  galerkin methods using finite volume sub-cells.
\newblock {\em Journal of Scientific Computing}, 70(3):1262--1289, 2017.

\bibitem{vilar2019posteriori}
Fran{\c{c}}ois Vilar.
\newblock A posteriori correction of high-order discontinuous galerkin scheme
  through subcell finite volume formulation and flux reconstruction.
\newblock {\em Journal of Computational Physics}, 387:245--279, 2019.

\bibitem{neumann50shocks}
John Von~Neumann and Robert Richtmyer.
\newblock A method for the numerical calculation of hydrodynamic shocks.
\newblock {\em Journal of Applied Physics}, 21:232--237, 1950.

\bibitem{colella84shocks}
Paul Woodward and Phillip Colella.
\newblock The numerical simulation of two-dimensional fluid flow with strong
  shocks.
\newblock {\em J. Comput. Phys.}, 54(1):115--173, 1984.

\end{thebibliography}

\section{Appendix}
\subsection{Proof of Lemma \ref{lem:solv}}
We first introduce some appropriate notation. Consider the set of multi-indices of total degree $p$ in $d=1,2,3$ dimensions given by
\[
	\N_p := \{ \bal \in \N_0^d \;| \;  |\bal| \le p\},
\]
where $|\bal| = \sum_{i=1}^d \alpha_i$ and define
\[
	\bx^\bal := \prod_{i=1}^d x_i^{\alpha_i},
	\qquad
	{{\bal}\choose{\bbe}} 
	:= \prod_{i=1}^d 
	{{\alpha_i}\choose{\beta_i}}  
\]
for any monomial $\bx=(x_1,\ldots,x_d)\in \R^d$. Here, we assume that ${{n}\choose{k}} =0$ for any $k>n$ so that ${{\bal}\choose{\bbe}} =0$ if there exists $i$ such that $\alpha_i<\beta_i$.
The corresponding space of $d$-variate polynomials of total degree at most $p$ given by
\[
	\mathbb P_{\! p} := \mbox{span}\{ \bx^\bal \; | \; \bal\in \N_p \},
\]
is of dimension $\mbox{dim}(\mathbb P_{\! p}) = {{p+d}\choose{d}}$ where the set of monomials $\bx^\bal$ consists of a basis.

We consider the case $r=p$ since the result follows immediately from the case $r=p$ for $r>p$.
Without loss of generality, we scale the unit simplex by a factor $p+1$ and therefore introduce
\[
	S_p = \{ \bx\in\R^d \; | \; \forall i=1,\ldots,d:  x_i\in[0,p+1] \mbox{ and } |\bx|\le p+1\}.
\]
Then, the sub-grid is given by the uniform partition of $S_p$ into $(p+1)^d$ simplices that are congurent to the unit simplex $\hat S := S_0$.
Indeed, $d!$ unit simplices fit into the unit cube in $d$ dimensions. Further, the cube $[0,p+1]^d$ consists of $(p+1)^d$ unit cubes which can be divided into $d!(p+1)^d$ unit simplices. Among all these unit simplices, there are $(p+1)^d$ unit simplices contained in $S_p$. 

We consider the subset of simplices 
\[
	S_\bal = \{ \hat\bx +\bal \; | \; \hat\bx \in  \hat S \}
\]
generated by different $\bal\in \N_p$. There are therefore also  $|\N_p| = \mbox{dim}(\mathbb P_{\! p})$ simplicies and they all belong to $S_p$.
We then consider the matrix 
\[
	A_{\bal\bbe}
	= \int_{S_\bal} \bx^\bbe \, d\bx
	= \int_{\hat S} (\hat \bx+\bal)^\bbe \, d\hat\bx
	= \sum_{\bga\in\N_p} {{\bbe}\choose{\bga}} \bal^{\bga}  \int_{\hat S} \hat \bx^{\bbe-\bga} \, d\hat\bx
	= \sum_{\bga\in\N_p} {{\bbe}\choose{\bga}} \bal^{\bga}  \bm f({\bbe-\bga})
\]
for all $\bal,\bbe\in\N_p$ and where $\bm f({\bal}) =  \int_{\hat S} \hat \bx^{\bal} \, d\hat\bx$. 
Here we have employed a change of variable from $S_\bal$ to $\hat S$ and the binomial formula. 
Now, introduce the two additional matrices
\[
	M_{\bal\bbe} = \bal^\bbe,
	\qquad\mbox{and}\qquad
	B_{\bal\bbe} = {{\bbe}\choose{\bal}}\bm f({\bbe- \bal}),
\]
so that $A=MB$. The matrix $M$ is indeed the interpolation matrix of the monomials $\bx^\bbe$ on the grid $\bal\in\N_p$, which is invertible. 
Second, the matrix $B$ is upper triangular with constant and non-zero diagonal 
\[
	B_{\bal\bal} = {{\bal}\choose{\bal}}\bm f({\bz}) = |\hat S|,
\]
and thus also invertible. In turn, $A$ is invertible. 

Let us now conclude the proof. 
Indeed, we consider any polynomial $v_p\in\mathbb P_{\! p}$ given by $v_p(\bx)=\sum_{\bbe \in \N_p } v_\bbe \,\bx^\bbe$ with $v_\bbe \in\mathbb R$. We now assume that $v_p$ has zero average on each sub-cell. Among those $n=(p+1)^d$ sub-cells, we only consider the set $S_\bal$ with $\bal\in\N_p$ and the conditions become: for all $\bal\in\N_p$, there holds
\[
	\int_{S_\bal} v_p(\bx) \, d\bx = 0 
	\qquad\Leftrightarrow\qquad
	\sum_{\bbe \in \N_p } v_\bbe \int_{S_\bal} \bx^\bbe\, d\bx = \sum_{\bbe \in \N_p } A_{\bal\bbe} v_\bbe = 0.
\]
Since $A$ is invertible, it follows that $v_\bbe=0$ for all $\bbe\in\N_p$ and the mapping $\ProjnK: \hVpK \to  \VnK$ is injective.
\end{document}